\documentclass{amsart}

\usepackage{amsmath}
\usepackage{amsfonts}
\usepackage{amssymb,enumerate}
\usepackage{amsthm}
\usepackage[all]{xy}

\newtheorem{lem}{Lemma}[section]
\newtheorem{cor}[lem]{Corollary}
\newtheorem{prop}[lem]{Proposition}
\newtheorem{thm}[lem]{Theorem}

\newtheorem{intthm}{Theorem}

\newtheorem{Defn}[lem]{Definition}
\newtheorem{Ex}[lem]{Example}
\newtheorem{Quest}[lem]{Question}
\newtheorem{Property}[lem]{Property}
\newtheorem{Properties}[lem]{Properties}
\newtheorem{Subprops}{}[lem]
\newtheorem{Para}[lem]{}
\newtheorem{Obs}[lem]{Observation}

\newtheorem{rem}[lem]{Remark}

\newenvironment{rmk}{\begin{rem}\rm}{\end{rem}}
\newenvironment{defn}{\begin{Defn}\rm}{\end{Defn}}
\newenvironment{ex}{\begin{Ex}\rm}{\end{Ex}}

\newenvironment{para}{\begin{Para}\rm}{\end{Para}}
\newenvironment{obs}{\begin{Obs}\rm}{\end{Obs}}

\theoremstyle{definition}

\newcommand{\ideal}[1]{\mathfrak{#1}}

\newcommand{\p}{\ideal{p}}

\newcommand{\ann}{\operatorname{Ann}}

\newcommand{\Hom}{\operatorname{Hom}}

\newcommand{\spec}{\operatorname{Spec}}
\newcommand{\A}{\mathcal{A}}
\newcommand{\X}{\mathcal{X}}
\newcommand{\B}{\mathcal{B}}
\newcommand{\I}{\mathcal{I}}
\newcommand{\F}{\mathcal{F}}
\newcommand{\PP}{\mathcal{P}}
\newcommand{\Rp}{R_{\p}}
\newcommand{\tor}{\operatorname{Tor}}
\newcommand{\ext}{\operatorname{Ext}}
\newcommand{\mcX}{\mathcal{X}}

\newcommand{\xra}{\xrightarrow}

\newcommand{\pd}{\operatorname{pd}}
\newcommand{\xpd}{\mcX\text{-}\pd}
\newcommand{\id}{\operatorname{id}}
\newcommand{\xid}{\mcX\text{-}\id}

\newcommand{\Ext}{\mathrm{Ext}_R}
\newcommand{\Tor}{\mathrm{Tor}^R}

\newcommand{\prim}{\mathfrak{p}}

\begin{document}

\author{Henrik Holm \ \ }

\address{ {\flushleft Department} of Mathematical Sciences, University
  of Aarhus, Ny Mun\-kegade, Building 1530,
  DK--8000 Aarhus C, Denmark}
\email{holm@imf.au.dk} \urladdr{http://home.imf.au.dk/holm/}

\author{\ \ Diana White}

\address{ {\flushleft Department} of Mathematics, University of Nebraska,
   203 Avery Hall, Lincoln, NE, 68588-0130 USA} \email{dwhite@math.unl.edu}
\urladdr{http://www.math.unl.edu/~s-dwhite14/}

\title[Foxby equivalence]{Foxby equivalence over
associative rings}





\begin{abstract}
  We extend the definition of a semidualizing module to associative
  rings.  This enables us to define and study Auslander and Bass classes with
  respect to a semidualizing bimodule $C$.  We then study the classes of
  $C$-flats, $C$-projectives, and $C$-injectives, and use them to
  provide a characterization of the modules in the Auslander and Bass
  classes. We extend Foxby equivalence to this new setting.  This
  paper contains a few results which are new in the commutative,
  noetherian setting.
\end{abstract}

\maketitle

\section*{Introduction}\label{sec:intro}

Over a commutative, noetherian local ring, semidualizing modules
provide a common generalization of a dualizing (canonical) module and
a free module of rank one. Foxby~\cite{foxby:gmarm} first defined them
(PG-modules of rank one), while Golod~\cite{golod:gdgpi} (suitable
modules) and Vasconcelos~\cite{vasconcelos:dtmc} (spherical modules)
furthered their study.  Recently,
Araya-Takahashi-Yoshino~\cite{araya:hiasdb} extended this definition
to a pair of non-commutative, but noetherian rings, while
White~\cite{white:gcproj} extended the definition to the
non-noetherian, but commutative, setting.  In this paper (see
Section~\ref{sec:sdm}), we define and study semidualizing
$(S,R)$-bimodules,
where $R$ and $S$ are arbitrary associative rings,
thereby encompassing all of the
aforementioned definitions.

Over a commutative noetherian ring, Avramov and
Foxby~\cite{avramov:rhfgd,foxby:gmarm} and
Enochs-Jenda-Xu~\cite{enochs:fdgipm} connected the study of
(semi)dualizing modules to associated Auslander and Bass classes for
(semi)dualizing modules, $\A_C(R)$ and $\B_C(R)$, which are
subcategories of the category of $R$-modules. This paper furthers
this study, which in our setting involves an Auslander class
$\A_C(R)$ and a Bass class $\B_C(S)$.

Those familiar with this area may wonder why we do not deal with
derived Auslander and Bass classes in this paper.  The short answer
is that the non-commutative situation is more subtle. A longer
answer is contained in Remark~\ref{dualizing}.

Many results
for Auslander and Bass classes associated to a semidualizing module
over a commutative noetherian ring, carry over to an associative
ring in a straightforward manner. However, some complications do
arise. Thus, in Section~\ref{sec:separating}, we define
\emph{faithfully} semidualizing bimodules, see Definition~\ref{sep}.
Over a commutative ring, all semidualizing modules are faithfully
semidualizing; see Theorem~\ref{commsep}, but it is unknown to the
authors if this is true in the non-commutative setting. We provide
many examples of ones that are, see Proposition~\ref{induction} and
Example~\ref{flatalgebras}.

Section~\ref{sec:abdef} contains basic properties of the Auslander
and Bass classes.  For example, the Auslander class $\A_C(R)$
contains all flat $R$-modules, and the Bass class $\B_C(S)$ contains
all injective $S$-modules; see Lemma~\ref{flat}. Also, both are
closed under summands, products, coproducts, and filtered colimits;
see Proposition~\ref{absums}.

To further the study of the classes $\A_C(R)$ and $\B_C(S)$, in
Section~\ref{sec:cproj} we introduce the classes of modules
$\F_C(S)$, $\PP_C(S)$, and $\I_C(R)$, consisting of the
 $C$-flats, $C$-projectives, and $C$-injectives, respectively. When
$C$ is free of rank one, these are the flats, projectives, and
injectives, respectively. The following is Theorem~\ref{Foxby
duality}.

\begin{intthm}[{\bf Foxby equivalence}]\label{intro:Foxby duality}
Let ${}_SC_R$ be a semidualizing bimodule. There are equivalences of
categories
  \begin{displaymath}
    \xymatrix@C=20ex{\PP_R(R) \ar@{^(->}[d] \ar@<0.8ex>[r]_-{\sim}^-{C
        \otimes_R-} & \PP_C(S) \ar@{^(->}[d]
      \ar@<0.8ex>[l]^-{\Hom_S(C,-)} \\
      \F_R(R) \ar@{^(->}[d] \ar@<0.8ex>[r]_-{\sim}^-{C
        \otimes_R-} & \F_C(S) \ar@{^(->}[d]
      \ar@<0.8ex>[l]^-{\Hom_S(C,-)} \\
      \A_C(R) \ar@<0.8ex>[r]_-{\sim}^-{C \otimes_R-} & \B_C(S)
      \ar@<0.8ex>[l]^-{\Hom_S(C,-)} \\
      \I_C(R) \ar@{^(->}[u] \ar@<0.8ex>[r]_-{\sim}^-{C \otimes_R-} &
      \I_S(S). \ar@<0.8ex>[l]^-{\Hom_S(C,-)} \ar@{^(->}[u]
    }
  \end{displaymath}
\end{intthm}
Propositions~\ref{closed} and~\ref{cprojext}, and Lemma~\ref{lemce}
discuss the closure properties of these classes.
Proposition~\ref{ce} pertains to their (pre)covering and
(pre)enveloping properties and includes some results that are new in
the commutative, noetherian setting. For example, when $R$ is
commutative noetherian the $C$-flats are preenveloping and the
$C$-injectives are precovering.

Section~\ref{classifyAB} contains two main results of the paper,
which provide alternative characterizations of the modules in the
Auslander and Bass classes in terms of the $C$-injectives,
$C$-projectives, and $C$-flats. Here is one.

\begin{intthm}
  \label{intro:acrchar}
  Let ${}_SC_R$ be a semidualizing bimodule.  An $R$-module $M$
  belongs to $\A_C(R)$ if and only if there exists a complex of $R$-modules
  \begin{eqnarray*}
    X= \ \cdots \longrightarrow P_1 \longrightarrow P_0
    \longrightarrow U^0 \longrightarrow U^1 \longrightarrow
    \cdots
  \end{eqnarray*}
  that satisfies the following conditions
  \begin{enumerate}[\quad\rm(a)]
  \item The complex $X$ is exact;
  \item Each $P_i$ is $R$-projective (or $R$-flat);
  \item Each $U^i$ is $C$-injective;
  \item There is an isomorphism $M \cong
    \operatorname{Coker}(P_1 \longrightarrow P_0)$; and
  \item The complex $C\otimes_R X$ is exact.
  \end{enumerate}
  Moreover, if \mbox{$M\in\A_C(R)$} then any complex constructed by
  splicing together an augmented projective (or flat) resolution of $M$ and an
  augmented proper $C$-injective coresolution of $M$ will satisfy the
  above properties.
\end{intthm}

When $C$ is faithfully semidualizing, there is a clear analogy with
the local setting. First, the Auslander class contains the modules
of finite flat dimension, while the Bass class contains the modules
of finite injective dimension. Second, if any two modules in a short
exact sequence are in the Auslander (respectively, Bass) class, then
so is the third; see Corollary~\ref{twothree}.


\section{Background} \label{sec:back}

This section introduces a number of assumptions, definitions,
notions, and results which will be used throughout the paper.

\begin{para}
  \label{setup}
  Throughout this work, $R$ and $S$ are fixed associative rings
  with unities and all $R$- or $S$-modules are understood to be
  \emph{unital left} $R$- or $S$-modules.
  Right $R$- or $S$-modules are identified with left modules over
  the \emph{opposite rings} $R^\mathrm{op}$ or $S^\mathrm{op}$.

  The suggestive notation ${}_SM_R$ is used to denote that $M$ is an
  \emph{$(S,R)$-bimodule}. This means that $M$ is both a left $S$- and
  a right $R$-module, and that these structures are compatible; namely
  \begin{displaymath}
    s(xr)=(sx)r \ \ \text{ for all } \ \ s\in S,\ r\in R, \text{ and }
    x\in M.
  \end{displaymath}
  Finally, if $\mcX$ is a class of, say, $R$-modules, then $\mcX^f$ is
  the subclass of all finitely generated $R$-modules in $\mcX$.
  Throughout this background section, $\mcX$ will denote a fixed class
  of $R$-modules.
\end{para}


\begin{para}
  \label{cxs}
  An \emph{$R$-complex} is a sequence of $R$-module homomorphisms
  \begin{displaymath}
    X = \cdots\xra{\partial^X_{n+1}}X_n\xra{\partial^X_n}
    X_{n-1}\xra{\partial^X_{n-1}}\cdots
  \end{displaymath}
  such that $\partial^X_{n-1}\partial^X_{n}=0$ for each integer $n$.
\end{para}

In this paper, all resolutions will be built from precovers, which
we now discuss.

\begin{para} \label{precover} Let $M$ be an $R$-module.  A
  homomorphism $\phi\colon X \to M$ with $X\in\X$ is an
  $\X$-\emph{precover} of $M$ if for every homomorphism $\psi\colon Y
  \to M$ with $Y\in\X$, there exists a homomorphism $f\colon Y\to X$
  such that $\phi f=\psi$.  If every $R$-module admits an
  $\X$-precover, then we say that the class $\X$ is
  \emph{precovering}.

  An $\X$-\emph{cover} of $M$ is an $\X$-precover $\phi\colon X \to M$
  with the additional property that any endomorphism $f\colon X \to X$
  with $\phi=\phi f$ must be an automorphism.  If every $R$-module
  admits an $\X$-cover, then we say that the class $\X$ is \emph{covering}.

  \emph{Preenvelopes} and \emph{envelopes} are defined dually;
   see~\cite{enochs:rha} for further details.
\end{para}

\begin{para}\label{resolutions}

If the class $\X$ is precovering, then for any $R$-module $M$, there
exists an \emph{augmented proper $\X$-resolution} of $M$, that is, a
complex
\begin{displaymath}
    X^+ = \cdots\xra{\partial^X_{2}}X_1 \xra{\partial^X_{1}}X_0\to
    M\to 0
  \end{displaymath}
such that $\Hom_R(Y,X^+)$ is exact for all $Y\in\X$. The truncated
complex
\begin{displaymath}
    X = \cdots\xra{\partial^X_{2}}X_1 \xra{\partial^X_{1}}X_0\to 0
  \end{displaymath}
is a \emph{proper $\X$-resolution} of $M$.

Note that $X^+$ need not be exact.  However, if $\X$ contains the
projectives, then any augmented proper $\X$-resolution of $M$ is
exact. \emph{Augmented proper $\X$-coresolutions} are defined
dually, and they must
  be exact if the class $\X$ contains the injectives.
\end{para}

\begin{para}\label{fpinf}
  A \emph{degreewise finite projective resolution} of an $R$-module
  $M$ is a projective resolution $P$ of $M$ such
  that each $P_i$ is finitely generated (projective).
\end{para}

\begin{para}
  \label{xdim}
  When $X$ is a precovering class, the \emph{$\mcX$-projective dimension} of $M$
  is
  \begin{displaymath}
    \xpd_R(M)= \inf \left\{ \sup\{n\mid X_n\neq 0\} \left|
      \begin{array}{l}
        \text{$X$ is an $\mcX$-resolution of $M$}
      \end{array}
     \right. \!\!\!
     \right\}
  \end{displaymath}
  The modules of $\X$-projective dimension zero are the non-zero
  modules in $\mcX$.  When $\X$ is preenveloping, the
  \emph{$\X$-injective dimension},
  denoted $\xid(-)$, is defined dually.
\end{para}

\begin{para}
  \label{extensions}
  Let $0 \to M' \to M \to M'' \to 0$ be an exact sequence of
  $R$-modules.  The class $\X$ is \emph{closed under extensions} if it
  has the property that if $M'$ and $M''$ are in $\X$, then so is $M$.
    The class $\X$ is
  \emph{closed under kernels of epimorphisms} if whenever $M$ and $M''$ are in $\X$,
  then so is $M'$.  Finally, the class $\X$ is \emph{closed under
  cokernels of monomorphisms} if whenever $M'$ and $M$ are in $\X$,
  then so is $M''$.
\end{para}

\begin{para}  \label{projres}
  The class $\mcX$ is \emph{projectively resolving} if
  \begin{enumerate}[\quad(a)]
  \item The class $\mcX$ contains every projective $R$-module, and
  \item For every exact sequence of $R$-modules $0\to M'\to M\to
    M''\to 0$ with $M''\in\mcX$, one has $M\in\mcX$ if and only if
    $M'\in\mcX$.
  \end{enumerate}
  The notion of \emph{injectively resolving} is defined dually.

  The class $\mcX$ is \emph{finite projectively resolving} if
  \begin{enumerate}[\quad(a)]
  \item The class $\mcX$ consists entirely of finitely generated $R$-modules,
  \item The class $\mcX$ contains every finitely generated projective
    $R$-module, and
  \item For every exact sequence of finitely generated $R$-modules
    $0\to M'\to M\to M''\to 0$ with $M''\in\mcX$, one has $M\in\mcX$
    if and only if $M'\in\mcX$.
  \end{enumerate}
\end{para}

The next result will be used in the proof of Theorem~\ref{commsep}.

\begin{para}\label{gruson}
  Vasconcelos \cite[(4.3)]{vasconcelos:dtmc} proves the following
  corollary of a theorem of Gruson: Let $R$ be a commutative ring and
  let $C$ be a finitely generated
  $R$-module with $\ann_R(C)=0$. If $M$ is any module such that $C\otimes_R M=0$,
  then $M=0$.
\end{para}

\begin{rmk}  \label{torext}
  Let ${}_SC_R$ be an $(S,R)$-bimodule.  Unless otherwise mentioned,
  an Ext group of the form $\ext^*_S({}_SM,C)$ will be computed
  by resolving ${}_SM$ with a projective resolution. This has the
  important consequence that $\ext^*_S({}_SM,C)$ inherits a
  right $R$-structure. Similar remarks could be said about the
  computation of other derived functors such as
  $\ext_{R^\mathrm{op}}^*(M_R,C)$ and
  $\tor_R^*(C,{}_RM)$.
\end{rmk}

We conclude this section with some necessary results about two
important homomorphisms. The proofs are straightforward, keeping
Remark~\ref{torext} in mind.

\begin{para}
  \label{tensorevaldefn}
  Consider modules ${}_SM$, ${}_SN_R$, and ${}_RF$.  The \emph{tensor
    evaluation} homomorphism
  \begin{displaymath}
    \omega_{MNP}\colon\Hom_S(M,N)\otimes_R F
    \longrightarrow \Hom_S(M,N\otimes_R F)
  \end{displaymath}
  is defined by $\omega_{MNF}(\psi\otimes f)(m)=\psi(m)\otimes f$.
  It is straightforward to verify that this is an isomorphism when
  $M$ is a finitely generated projective.
  In general, $\omega_{MNF}$ is just an abelian group homomorphism.
  However, if ${}_SM$ has an additional right $R$-structure
  compatible with the given left $S$-structure, then $\omega_{MNF}$
  becomes $R$-linear.
\end{para}

\begin{lem}  \label{tensoreval}
  Let ${}_SM$, ${}_SN_R$, and ${}_RF$ be modules such that $M$
  admits a degreewise finite $S$--projective resolution, and let $F$
  be flat. Let $i$ be an integer.
  \begin{enumerate}[\quad\rm(a)]
  \item\label{t2} The map $\omega_{MNF}$ induces an isomorphism of abelian
    groups
    \begin{displaymath}
      \ext^i_S(M,N)\otimes_R F \cong
      \ext^i_S(M,N\otimes_RF).
    \end{displaymath}
  \item\label{t3} If $\,\ext^i_S(M,N)=0$, then
    $\ext^i_S(M,N\otimes_R F)=0$.
  \item\label{t4} If $F$ is faithfully flat and
    $\ext^i_S(M,N\otimes_R F)=0$, then
    \mbox{$\ext^i_S(M,N)=0$.}\qed
  \end{enumerate}
\end{lem}

\begin{para}
  \label{homevaldefn}
  Consider modules $M_R$, ${}_SN_R$, and ${}_SI$.  The
  \emph{Hom-evaluation} homomorphism
  \begin{displaymath}
    \theta_{MNI}\colon M_R \otimes_R\Hom_S(N,I) \longrightarrow
    \Hom_S(\Hom_{R^\mathrm{op}}(M,N),I)
  \end{displaymath}
  is defined by $\theta_{MNI}(m\otimes
  \phi)(\varphi)=(\phi\circ\varphi)(m)$.  It is straightforward to
  verify that this is an isomorphism when $M$ is a finitely generated
  projective.
  In general, $\theta_{MNI}$ is just a homomorphism of abelian groups,
  but if $M_R$ has an additional left $S$-structure
  compatible with the right $R$-structure, then $\theta_{MNI}$
  becomes $S$-linear.
\end{para}

\begin{lem}
  \label{homeval}
  Let $M_R$, ${}_SN_R$, and ${}_SI$ be modules such that $M$ admits a
  degreewise finite $R^\mathrm{op}$--projective resolution, and let $I$ be
  injective. Let $i$ be an integer.
  \begin{enumerate}[\quad\rm(a)]
  \item\label{h2} The map $\theta_{MCI}$ induces an isomorphism of abelian
    groups
    \begin{displaymath}
      \tor^R_i(M,\Hom_S(N,I)) \cong
      \Hom_S(\ext_{R^\mathrm{op}}^i(M,N),I).
    \end{displaymath}
  \item\label{h3} If $\,\ext^i_{R^\mathrm{op}}(M,N)=0$, then
    \mbox{$\tor^R_i(M,\Hom_S(N,I))=0$}.
  \item\label{h4} If $I$ is faithfully injective and
    \mbox{$\tor^R_i(M,\Hom_S(N,I))=0$}, then there is an equality
    $\ext_{R^\mathrm{op}}^i(M,N)=0$.\qed
  \end{enumerate}
\end{lem}

\section{Semidualizing bimodules}\label{sec:sdm}

We begin by extending the definition of a semidualizing module.

\begin{defn}  \label{sdm}
  An $(S,R)$-bimodule $C={}_SC_R$ is \emph{semidualizing} if
  \begin{enumerate}[\quad\rm(a)]
  \item[(a1)] ${}_SC$ admits a degreewise finite $S$-projective
    resolution.
  \item[(a2)] $C_R$ admits a degreewise finite
    $R^\mathrm{op}$-projective resolution.
  \item[(b1)] The homothety map ${}^{}_SS_S
  \xra{{}_S\gamma}
    \Hom_{R^\mathrm{op}}(C,C)$ is an isomorphism.
  \item[(b2)] The homothety map ${}^{}_RR_R
  \xra{\gamma_R}
    \Hom_S(C,C)$ is an isomorphism.
  \item[(c1)] $\ext^{\geqslant 1}_S(C,C)=0$.
  \item[(c2)] $\ext^{\geqslant
      1}_{R^\mathrm{op}}(C,C)=0$.
  \end{enumerate}
\end{defn}


\begin{para}\label{sym}
Unless otherwise stated, when $R=S$ is commutative, all
semidualizing bimodules in this paper are \textsl{symmetric} in the
sense that the two $R$-actions on $C$ agree.  In this case we will
use the terminology ``$C$ is semidualizing over $R$''.  Note that
when $R=S$ is commutative and noetherian, Definition~\ref{sdm}
agrees with the established terminology; that is a finitely
generated
  $R$-module $C$ is \emph{semidualizing} if the natural homothety map
  $R\xra{\gamma_R}\Hom_R(C,C)$ is an isomorphism and $\ext^{\geqslant
    1}_R(C,C)=0$.  Two examples are the free module of rank 1, and over a
    Cohen-Macaulay local ring, the dualizing (canonical)
  module, when it exists.
\end{para}

\begin{obs}\label{sdmobs}
  \begin{enumerate}[(a)]
  \item When $S$ is left noetherian and $R$ is right noetherian,
    conditions (a1) and (a2) reduce to ${}_SC$ and $C_R$ being
    finitely generated, and therefore Definition~\ref{sdm} agrees with
    that of Araya-Takahashi-Yoshino~\cite[(2.1)]{araya:hiasdb}.
  \item Let $R=S$ be commutative.  Conditions (a1) and (a2) reduce to
    the statement that $C$ admits a degreewise finite projective
    resolution, while conditions (b1) and (b2) reduce to
    $\Hom_R(C,C)\cong R$, and conditions (c1) and (c2) reduce to
    $\ext_R^{\geqslant 1}(C,C)=0$.  Thus, Definition~\ref{sdm} agrees
    with that of White~\cite[(1.8)]{white:gcproj}.
  \item By the symmetry of the definition above it is clear that if
    $C$ is a semidualizing $(S,R)$-bimodule, then $C$ is also a
    semidualizing $(R^\mathrm{op},S^\mathrm{op})$-bimodule.
  \item \label{sdmobs4} If $R=S$ is commutative and $C$ is a
    semidualizing $R$-module, then $C_{\p}$ is a semidualizing
    $\Rp$-bimodule for all $\p\in\spec(R)$, cf.~Example
    \ref{Morita}(a) below.
\end{enumerate}
\end{obs}

Note that~\ref{Morita}(a) below can be applied to produce examples
of semidualizing modules over a commutative, but not necessarily
noetherian, ring. Many additional examples can be found in the
next section.

\begin{ex}  \label{Morita}
  \begin{enumerate}[(a)]
  \item Let $Q\longrightarrow R$ be a flat ring homomorphism between
    commutative rings.  If $E$ is semidualizing over $Q$, then
    $E\otimes_Q R$ is semidualizing over $R$, cf.~Proposition
    \ref{induction}.
  \item Assume the $R$ and $S$ are Morita equivalent rings in the
    sense of \cite[(9.5)]{weibel:iha}, that is, there exists bimodules
    ${}_RP_S$ and ${}_SQ_R$ with
  \begin{displaymath}
    {}_RP_S \otimes_S {}_SQ_R \cong {}_RR_R
    \ \text{ and } \
    {}_SQ_R \otimes_R {}_RP_S \cong {}_SS_S.
  \end{displaymath}
  By~\cite[(9.5.4), (9.5.2)]{weibel:iha}, the bimodule
  ${}_SQ_R$ is semidualizing.
\end{enumerate}
\end{ex}

We close this section with a comparison of \emph{derived Auslander and
  Bass classes} \cite{christensen:sdctac} and \emph{module Auslander
  and Bass classes} \cite{foxby:gmarm} in the non-commutative setting.

\begin{rmk}
  \label{dualizing}
  Over a commutative, noetherian ring $R$,
  Christensen~\cite{christensen:sdctac} generalized the notions of a
  semidualizing module and a dualizing complex to that of a
  semidualizing complex $C$.  He
  then connected their study to the associated \emph{derived} Auslander
    and Bass classes, $\mathsf{A}_C(R)$ and $\mathsf{B}_C(R)$, which
  are subcategories of the derived category of $R$.
  When $C$ is a semidualizing module, Foxby \cite{foxby:gmarm}
  studied the \emph{module} Auslander and Bass classes, $\A_C(R)$ and
  $\B_C(R)$, which are subcategories of the category of $R$-modules
  (see Definition \ref{aus}).

  Thus, for a semidualizing module $C$, a natural question arises:
  Does $\A_C(R)$ and $\B_C(R)$ consist of the modules belonging to
  $\mathsf{A}_C(R)$ and $\mathsf{B}_C(R)$?

  In the commutative noetherian setting, the answer is ``yes'', and
  thus the existing literature tends to focus on the more general
  derived Auslander and Bass classes. However, for non-commutative
  rings the question is much more complicated:

  Dualizing complexes of bimodules over a non-commutative but two-sided
  noetherian ring have been given several different definitions,
  e.g.~\cite{frankild:ddgmgdga,miyachi:dcmdt,yekutieli:radc}.
  In~\cite{CFH} the authors use a variant of Miyachi's
  definition~\cite{miyachi:dcmdt} of a dualizing complex of bimodules to
  consider a special case of the derived Auslander and Bass classes.
  However, this definition requires the existence of a so-called
  \emph{biprojective resolution} of the dualizing complex of bimodules, and
  such a resolution is only known to exist in certain special cases.
  Even if the semidualizing bimodule admits a
  biprojective resolution, so that the derived Auslander and Bass
  classes can be defined, it is not known to these authors if the
  modules in $\mathsf{A}_C(R)$ or $\mathsf{B}_C(R)$ belong to
  $\A_C(R)$ or $\B_C(R)$?
\end{rmk}

\section{Faithfully semidualizing bimodules}\label{sec:separating}

 This section focuses on \emph{faithfully} semidualizing
 bimodules.  Over a commutative ring, all semidualizing modules are
 faithfully semidualizing; see Theorem~\ref{commsep}.

\begin{defn}  \label{sep}
  A semidualizing bimodule ${}_SC_R$ is \emph{faithfully semidualizing} if it
  satisfies the following conditions for all modules ${}_SN$ and $M_R$.
  \begin{enumerate}[\quad\rm(a)]
  \item If $\Hom_S(C,N)=0$, then $N=0$.
  \item If $\Hom_{R^\mathrm{op}}(C,M)=0$, then $M=0$.
  \end{enumerate}
\end{defn}

\begin{para}
  \label{sepcommdefn}
  Note that when $R$ is commutative, these conditions are the same,
  and a semi\-dualizing module $C$ is faithfully semidualizing if and
  only if, for any $R$-module $N$, the condition $\Hom_R(C,N)=0$
  implies $N=0$.
\end{para}

\begin{para}
  By left exactness of $\Hom_S(C,-)$ and $\Hom_{R^\mathrm{op}}(C,-)$,
  it suffices to verify that the conditions in Definition \ref{sep} hold for
  all \emph{cyclic} modules $N$ and $M$.
\end{para}

\begin{lem}
  \label{defininglemma}
  A semidualizing bimodule ${}_SC_R$ is faithfully semidualizing if
  and only if the following conditions hold for all modules $N_S$ and
  ${}_RM$.
    \begin{enumerate}[\quad\rm(a)]
    \item If $N \otimes_S C=0$, then $N=0$.
    \item If $C \otimes_R M=0$, then $M=0$.
    \end{enumerate}
\end{lem}

\begin{proof}
  Assume that $N \otimes_SC=0$. This provides the equality
  below
  \begin{displaymath}
    \Hom_S(C,\Hom_{\mathbb{Z}}(N,\mathbb{Q}/\mathbb{Z}))
    \cong
    \Hom_{\mathbb{Z}}(N \otimes_SC,\mathbb{Q}/\mathbb{Z})=0
  \end{displaymath}
  while the isomorphism follows from Hom-tensor adjointness.  Since
  $C$ is faithfully semidualizing,
  $\Hom_{\mathbb{Z}}(N,\mathbb{Q}/\mathbb{Z})=0$.  The module
  $\mathbb{Q}/\mathbb{Z}$ is faithfully injective, which forces
  $N=0$.  Similarly, if $C\otimes_RM=0$, then
  $M=0$.

  The converse is proved similarly, using
  Lemma~\ref{homeval}\eqref{h2}. 
\end{proof}

\enlargethispage{.9cm}

\begin{para}
  When $R$ is commutative and noetherian, it is straightforward to
  prove that every semidualizing module $C$ is faithfully
  semidualizing.  Indeed, if $M\neq 0$, then for any associated prime
  $\p$ of $M$ there is a nonzero map $C_{\p}\to k(\p) \to M_{\p}$.  It
  follows that there is a nonzero map $C\to M$.
\end{para}

The situation where $R$ is commutative, but non-noetherian seems to
require significantly more work (which is done behind the scenes in
Gruson's Theorem; see~\ref{gruson}), but the same result holds, as
we show next. The authors are unaware if this result carries over to
the non-commutative setting.  However, Proposition~\ref{induction}
and Example~\ref{flatalgebras} below provide many examples of
faithfully semidualizing bimodules over a wide class of
non-commutative rings.

\begin{prop}\label{commsep}
  Let $R$ be a commutative ring.  If $\,C$ is a semidualizing
  $R$-module, then $C$ is faithfully semidualizing.
\end{prop}

\begin{proof}
  As $R\cong\Hom_R(C,C)$ we have $\ann_R(C)=0$. Therefore,
  Lemma~\ref{defininglemma} and the corollary of Gruson's
  Theorem~\ref{gruson} imply that $C$ is faithfully semidualizing.
\end{proof}

The next result is a non-commutative, module version of
\cite[(5.1)]{christensen:sdctac}.

\begin{prop}  \label{induction}
  Let $Q$ be a commutative ring and let $R$ be an associative and flat
  $Q$-algebra. If $E$ is a  semidualizing
  $Q$-module, then \mbox{${}_RC_R=E\otimes_Q{}_RR_R$} is a
  faithfully semidualizing $(R,R)$-bimodule.
\end{prop}

\begin{proof}
  We first prove that $C$ is a semidualizing bimodule.  We show only
  that it satisfies (a1), (b1), and (c1) of Definition~\ref{sdm}, as
  the other parts are proved similarly.

  Since $E$ admits a degreewise finite $Q$-projective resolution $P$ and $R$ is
  $Q$-flat, $P\otimes_Q{}_RR$ is a degreewise finite projective
  resolution of ${}_RC=E\otimes_Q{}_RR$.

  The map \mbox{$Q \xra{\gamma_Q} \Hom_Q(E,E)$} is an isomorphism, so
  the commutative diagram
  \begin{displaymath}
    \xymatrix{
      R \ar[r]^-{{}^{}_R\gamma}
      \ar[d]^-{\cong}_-{\mathrm{}} &
      \Hom_{R^\mathrm{op}}(C,C) \ar@{=}[r]^-{} &
      \Hom_{R^\mathrm{op}}(E\otimes_QR,E\otimes_QR)
      \ar[d]^-{\textrm{adjointness}}_-{\cong}
       \\
      Q \otimes_Q R
      \ar[d]^-{\cong}_-{\gamma^{}_Q\otimes_QR} &
      {} &
      \Hom_Q(E,\Hom_{R^\mathrm{op}}(R,E\otimes_QR))
      \ar[d]_-{\cong}^-{}
       \\
      \Hom_Q(E,E)\otimes_QR
      \ar[rr]^-{\cong}_-{\omega_{EER}} &
      {} &
      \Hom_Q(E,E\otimes_QR)
    }
  \end{displaymath}
  shows that the homomorphism \mbox{$R
    \xra{{}^{}_R\gamma} \Hom_R(C,C)$} is an isomorphism.

  Finally, let $P$ be a $Q$-projective resolution of $E$.  For $i>0$, the
  first and fourth isomorphisms below are by definition while the second
   is Hom-tensor adjointness
  \begin{align*}
    \ext^i_R(C,C)
     &=
    \operatorname{H}_{-i}\Hom_R(P\otimes_QR,E\otimes_QR) \\
     &\cong
    \operatorname{H}_{-i}\Hom_Q(P,\Hom_R(R,E\otimes_QR)) \\
     &\cong
    \operatorname{H}_{-i}\Hom_Q(P,E\otimes_Q R) \\
     &\cong
    \ext^i_Q(E,E \otimes_Q R)=0.
  \end{align*}
  The third isomorphism is standard, and the last equality is by
  Lemma~\ref{tensoreval}\eqref{t3}, as $\ext^i_Q(E,E)=0$ and $R$ is
  $Q$-flat.  Next, note that for modules $N_R$ and ${}_RM$ there are
  standard isomorphisms
  \begin{displaymath}
    N\otimes_RC =
    N\otimes_R(E\otimes_QR) \cong
    (N\otimes_RR)\otimes_QE \cong N \otimes_Q
    E
  \end{displaymath}
  and similarly $C \otimes_R M \cong E \otimes_QM$.  Since $E$
  is faithfully semidualizing by Proposition \ref{commsep}, it follows from
  Lemma~\ref{defininglemma} that $C$ is faithfully semidualizing as well.
\end{proof}

\begin{ex} \label{flatalgebras} Here are several standard examples of
  associative flat algebras $R$ over a commutative ground ring $Q$.
  \begin{enumerate}[\rm(a)]

  \item If $T$ is any commutative and torsion-free (that is, flat)
    $\mathbb{Z}$-algebra, then the \emph{tensor product algebra}
    \begin{displaymath}
      R=T\otimes_{\mathbb{Z}}Q
    \end{displaymath}
    is a commutative flat $Q$-algebra, as $R \otimes_Q- \cong
    T\otimes_{\mathbb{Z}}-$. Even if $Q$ is noetherian, $R$ need not
    be.  For example, $T=\mathbb{Z}[x_1,x_2,x_3,\ldots]$ gives
    $R=Q[x_1,x_2,x_3,\ldots]$.

  \item The \emph{$n\times n$ matrix ring}
    \begin{displaymath}
      R=M_n(Q)
    \end{displaymath}
    is a free $Q$-algebra. In general, $Q$ is non-commutative.
    Moreover, $R$ is left and/or right noetherian if and only if $Q$
    is noetherian; see~\cite[(1.1.2)]{mcconnell:nnr}.

  \item Let \mbox{$\alpha\colon Q \longrightarrow Q$} be a ring
    endomorphism, and let \mbox{$\delta \colon Q \longrightarrow Q$}
    be an \emph{$\alpha$-derivation}. Then the \emph{skew polynomial
      ring}
    \begin{displaymath}
      R=Q[\theta;\alpha,\delta]
    \end{displaymath}
    is a free $Q$-algebra with basis $\{1,\theta,\theta^2,\ldots\}$.
    In general, $R$ is non-commutative. Even if $Q$ is noetherian, $R$
    need not be (from either side). However, if $Q$ is noetherian and
    $\alpha$ is an automorphism, then $R$ is two-sided noetherian;
    see~\cite[Exer. 1N and (1.12)]{goodearl:innr}.

  \item If $G$ is any group, then the \emph{group ring}
    \begin{displaymath}
      R=Q[G]
    \end{displaymath}
    is a free $Q$-algebra with basis $G$. Note that $R$ is commutative
    if and only if $G$ is abelian. Even if $Q$ is noetherian $R$ need
    not be (from either side). However, if $Q$ is noetherian and $G$
    is \emph{polycyclic-by-finite} then $R$ is two-sided noetherian;
    see~\cite[(1.5.12)]{mcconnell:nnr}.

  \item Let $F$ be any flat $Q$-module. Then the \emph{formal triangular
    matrix ring}
    \begin{displaymath}
      R=\left(
        \begin{array}{cc}
          Q & F \\
          0 & Q
        \end{array}
      \right)
    \end{displaymath}
    is a flat $Q$-algebra, since $R\cong Q\oplus F\oplus Q$ as a
    $Q$-module. In general, $R$ is non-commutative. Moreover $R$ is
    left and/or right noetherian if and only if $Q$ is noetherian and
    $F$ is finitely generated; see \cite[Exer. 1C and
    (1.9)]{goodearl:innr}.
  \end{enumerate}
\end{ex}

\section{Auslander and Bass classes with respect to $C$}\label{sec:abdef}

In this section we introduce and investigate properties of the
Auslander and Bass classes 
with respect to a semidualizing $(S,R)$-bimodule $C={}_SC_R$. Over a
commutative noetherian ring the following definition can be found in
\cite[sec.~1]{foxby:gmarm}.

\begin{defn}
  \label{aus}
  The \textit{Auslander class} $\A_C(R)$ with respect to $C$ consists
  of all $R$-modules $M$ satisfying
  \begin{enumerate}
  \item[(A1)] $\tor_{\geqslant 1}^R(C,M)=0$,
  \item[(A2)] $\ext^{\geqslant 1}_S(C,C\otimes_RM)=0$, and
  \item[(A3)] The natural evaluation homomorphism $\mu^{}_M\colon
    M\longrightarrow \Hom_S(C,C\otimes_R M)$ is an
    isomorphism (of $R$-modules).
  \end{enumerate}
   The \textit{Bass class} $\B_C(S)$ with respect to $C$
  consists of all $S$-modules $N$ satisfying
  \begin{enumerate}
  \item[(B1)] \label{bass1} $\ext^{\geqslant 1}_S(C,N)=0$,
  \item[(B2)] \label{bass2} $\tor_{\geqslant
      1}^R(C,\Hom_R(C,N))=0$, and
  \item[(B3)] \label{bass3} The natural evaluation homomorphism $\nu^{}_N\colon
    C\otimes_R\Hom_S(C,N)\longrightarrow N$ is an
    isomorphism (of $S$-modules).
  \end{enumerate}
\end{defn}

It is straightforward to check the following:

\begin{obs}
  \label{leftrightinverses}
  Given modules ${}_RM$ and ${}_SN$, the homomorphisms
  \begin{displaymath}
    \xymatrix{
      C\otimes_R\Hom_S(C,C\otimes_RM)
     \ar@<0.6ex>[rr]^-{\nu^{}_{(C\otimes_RM)}} & &
      C\otimes_RM
      \ar@<0.6ex>[ll]^-{C\otimes_R\mu^{}_M}, \text{ and}
    }
  \end{displaymath}
  \vspace{-2ex}
  \begin{displaymath}
    \xymatrix{
      \Hom_S(C,N) \ar@<0.6ex>[rr]^-{\mu^{}_{\Hom_S(C,N)}} & &
      \Hom_S(C,C\otimes_R\Hom_S(C,N))
      \ar@<0.6ex>[ll]^-{\Hom_S(C,\nu^{}_N)}
    }
  \end{displaymath}
from Definition~\ref{aus} yield identities
  \begin{eqnarray*}
    &\nu^{}_{(C\otimes_RM)} \circ (C\otimes_R\mu^{}_M) =
    \mathrm{id}_{(C\otimes_RM)} \text{ and }& \\
   &\Hom_S(C,\nu^{}_N) \circ \mu^{}_{\Hom_S(C,N)} =
    \mathrm{id}_{\Hom_S(C,N)}.&
  \end{eqnarray*}
\end{obs}

The next result is a  component of the Foxby equivalence; see
Theorem~\ref{Foxby duality}.

\begin{prop}  \label{abequiv}
  Let ${}_SC_R$ be a semidualizing bimodule. There are equivalences
  of categories
  \begin{displaymath}
    \xymatrix@C=20ex{
      \A_C(R) \ar@<0.8ex>[r]_-{\sim}^-{C \otimes_R-} & \B_C(S).
      \ar@<0.8ex>[l]^-{\Hom_S(C,-)}
    }
  \end{displaymath}
\end{prop}

\begin{proof}
  To see that the functor $C\otimes_R-$ maps $\A_C(R)$ to
  $\B_C(S)$, let \mbox{$M\in\A_C(R)$} and note that the
  following hold by definition of the class $\A_C(R)$.
  First, for all $i>0$ there is an equality
$\ext_S^i(C,C\otimes_RM)=0$.
  Second, for $i>0$ one has
  \begin{align*}
    0 =
    \tor_i^R(C,M)
     \cong
    \tor_i^R(C,\Hom_S(C,C\otimes_R M)).
  \end{align*}
  Finally, as $\mu^{}_M$ is an isomorphism, so is $C\otimes_R\mu^{}_M$,
  Thus, Observation~\ref{leftrightinverses} implies that
  $\nu^{}_{(C\otimes_RM)}$ is also an isomorphism and the inclusion
  $C\otimes_R M \in \B_C(S)$ follows.

  The proof that  $\Hom_S(C,-)$ maps $\B_C(S)$ to
  $\A_C(R)$ is similar.
  To conclude, note that if $M \in \A_C(R)$ and $N \in \B_C(S)$,
  then there are natural isomorphisms
  \begin{eqnarray*}
    &\mu^{}_M\colon M \stackrel{\cong}{\longrightarrow}
    \Hom_S(C,C\otimes_R M) \text{ and}& \\
    &\nu^{}_N\colon
    C\otimes_R\Hom_S(C,N)
    \stackrel{\cong}{\longrightarrow} N.
  \end{eqnarray*}
  The desired equivalences of categories now follows.
\end{proof}

The next lemma, together with Corollary~\ref{inclusions} (in view of
Proposition~\ref{commsep}), extend~\cite[(1.2)]{foxby:gmarm}.

\begin{lem}  \label{flat}
  Let ${}_SC_R$ be a semidualizing bimodule.
  The class $\A_C(R)$ contains the flat $R$-modules and the class
  $\B_C(R)$ contains the injective $S$-modules.
\end{lem}

\begin{proof}
  For a flat $R$-module $F$, one has
 $   \Tor_{\geqslant 1}(C,F)=0$ and
   Lemma~\ref{tensoreval}\eqref{t3} implies
  \begin{displaymath}
    \ext^{\geqslant 1}_S(C,C\otimes_R F)=0.
  \end{displaymath}
  Finally, Lemma~\ref{tensoreval}\eqref{t2} implies that
  $\omega_{CCF}$ is an isomorphism of abelian groups.  It is
  also $R$-linear, and so an $R$-module isomorphism.  Thus, the
  commutative diagram
  \begin{eqnarray*}
    \xymatrix{      R\otimes_R F \ar[d]_-{\cong}
      \ar[rr]_-{\cong}^-{\gamma_R \otimes_R\, F}
      & {} & \Hom_S(C,C)\otimes_R {}_RF
      \ar[d]^-{\omega_{CCF}}_-{\cong} \\
      F \ar[rr]_-{\mu^{}_F} & {} &
      \Hom_S(C,C\otimes_R F)
    }
  \end{eqnarray*}
  implies that $\mu^{}_F$ is an isomorphism.  The other proof is dual.
\end{proof}

We conclude this section by investigating how the Auslander and Bass
classes behave with respect to summands, products, coproducts, and
filtered colimits.

\begin{prop}  \label{absums}
  Let ${}_SC_R$ be a semidualizing bimodule.
  \begin{enumerate}[\quad\rm(a)]
  \item\label{absums1} The classes $\A_C(R)$ and $\B_C(S)$ are closed
    under direct summands, products, coproducts, and filtered colimits.
  \item\label{absums3} The classes $\A_C^f(R)$ and $\B_C^f(S)$ are
    closed under finite direct sums and direct summands.
  \end{enumerate}
\end{prop}

\begin{proof}
  \eqref{absums1} It is straightforward to verify that $\A_C(R)$ is
  closed under direct summands and finite direct sums, as the functors
  $\ext$, $\tor$, and $\Hom$ are additive. We prove that $\A_C(R)$ is
  closed under filtered colimits, from which it follows that $\A_C(R)$
  is closed under coproducts. To this end, let $\{M_\lambda\}$ be a
  filtered direct system of $R$-modules. Since tensor products and taking
  homology commute with filtered colimits, for each $i\geqslant 0$
  there is an isomorphism of abelian groups
  \begin{displaymath}
    \tag{\text{$\dagger$}}
    \tor_i^R(C,\varinjlim M_\lambda) \cong
    \varinjlim \tor_i^R(C,M_\lambda)
  \end{displaymath}
  and, when $i=0$, an isomorphism of $S$-modules
  \begin{displaymath}
    \tag{\text{$\dagger\!\dagger$}}
    C \otimes_R (\varinjlim M_\lambda) \cong
    \varinjlim (C \otimes_R M_\lambda).
  \end{displaymath}
  As ${}_SC$ admits a degreewise finite $S$-projective resolution,
  $\ext_S^i(C,-)$ commutes with filtered colimits.  In
  particular, for $i\geqslant 0$ the isomorphism in
  $(\dagger\!\dagger)$ gives an isomorphism of abelian groups
  \begin{displaymath}
    \tag{\text{$\dagger\!\!\dagger\!\!\dagger$}}
    \ext_S^i(C,C \otimes_R (\varinjlim M_\lambda)) \cong
    \varinjlim \ext_S^i(C,C \otimes_R M_\lambda)
  \end{displaymath}
  and, when $i=0$, an isomorphism of $R$-modules
  \begin{displaymath}
    \Hom_S(C,C \otimes_R (\varinjlim M_\lambda)) \cong
    \varinjlim \Hom_S(C,C \otimes_R M_\lambda).
  \end{displaymath}
  This isomorphism fits into the following commutative diagram
  \begin{displaymath}
    \xymatrix@C=-5ex{\varinjlim M_\lambda \ar[dr]_(0.4){\varinjlim
        (\mu^{}_{M_\lambda})} \ar[rr]^-{\mu^{}_{(\varinjlim
          M_\lambda)}} &  & \Hom_S(C,C \otimes_R
      (\varinjlim M_\lambda)) \\
    {} & \varinjlim \Hom_S(C,C \otimes_R M_\lambda)
    \ar[ur]_-{\cong} & {}
    }
  \end{displaymath}
  Since $\varinjlim (\mu^{}_{M_\lambda})$ is an isomorphism if each
  $\mu^{}_{M_\lambda}$ is, the diagram above, together with the
  isomorphisms $(\dagger)$ and $(\dagger\!\!\dagger\!\!\dagger)$,
  imply that $\A_C(R)$ is closed under filtered colimits.

A proof similar to the one above, using the proof
of~\cite[(3.2.26)]{enochs:rha}, shows that $\A_C(R)$ is closed under
products.
Similar arguments show that the Bass class $\B_C(S)$ is closed under
direct summands, coproducts, products, and filtered colimits.

  \eqref{absums3} Since finite direct sums and direct summands of
  finitely generated modules are finitely generated, this follows from~\eqref{absums1}.
\end{proof}

\section{$C$-flats, $C$-projectives, and $C$-injectives}\label{sec:cproj}

In this section we study flat, projective, and injective modules
with respect to a semidualizing bimodule $C={}_SC_R$, investigate
basic properties of these classes, and prove a form of Foxby
equivalence. In the commutative noetherian setting, the following
classes of modules already appear in, for
example,~\cite{enochs:fdgipm}, \cite{EEESY}, \cite{holm:sdmrghd},
and \cite{xu:fcm}.

\begin{defn}  \label{cproj}
  An $S$-module is \emph{$C$-flat} (resp.,
  \emph{$C$-projective}) if it has the form $C\otimes_RF$ for some
  flat (resp., projective) module ${}_RF$.  An $R$-module
  is \emph{$C$-injective} if it has the form $\Hom_S(C,I)$
  for some injective module ${}_SI$. Set the notation
  \begin{align*}
    &\F_C=\F_C(S)=\{C\otimes_RF \mid {}_RF \text{ is flat}\}, \\
    &\PP_C=\PP_C(S)=\{C\otimes_RP \mid {}_RP \text{ is
      projective}\}, \\
    &\I_C=\I_C(R)=\{\Hom_S(C,I) \mid {}_SI \text{ is
      injective}\}.
  \end{align*}
\end{defn}

It is straightforward to prove the following:

\begin{lem}\label{uw}Let  ${}_SC_R$ be a
  semidualizing bimodule.  For modules ${}_RU$ and ${}_SV$ the
  following hold.
  \begin{enumerate}[\quad\rm(a)]
  \item\label{FC} \mbox{$V \in \F_C(S) \iff V \in \B_C(S)$}
    and \mbox{$\Hom_S(C,V)$} is flat over $R$.
  \item\label{PC} \mbox{$V \in \PP_C(S) \iff V \in \B_C(S)$}
    and \mbox{$\Hom_S(C,V)$} is projective over $R$.
  \item\label{IC} \mbox{$U \in \I_C(R) \iff U \in \A_C(R)$}
    and \mbox{$C\otimes_RU$} is injective over $S$.
  \end{enumerate}
\end{lem}


The next result is a non-commutative, non-noetherian version
of~\cite[(1.4)]{foxby:gmarm}.
\begin{thm}[{\bf Foxby equivalence}]\label{Foxby duality}
Let ${}_SC_R$ be a semidualizing bimodule. There are equivalences of
categories as illustrated in Theorem~\ref{intro:Foxby duality} from
the introduction.
\end{thm}

\begin{proof}The equivalence between $\A_C(R)$ and $\B_C(S)$ was
established in Proposition~\ref{abequiv}.  The vertical containments
are either trivial or follow from Lemmas~\ref{flat} and~\ref{uw}.
The horizontal equivalences follow from Lemma~\ref{uw}.
\end{proof}

\begin{rmk}
  \label{FDsep}
  When $C$ is faithfully semidualizing, a stronger version of the
  above theorem holds, proved using Corollaries~\ref{inclusions2}
  and~\ref{inclusions}. Specifically, replace the classes $\PP_R(R)$,
  $\F_R(R)$ and $\I_S(S)$ with the classes of modules of finite
  $R$-projective, $R$-flat and $S$-injective dimension and replace the
  classes $\PP_C(S)$, $\F_C(S)$ and $\I_C(R)$ with the classes of
  modules of finite $\PP_C$-projective dimension over $S$, finite
  $\F_C$-projective dimension over $S$ and finite $\I_C$-injective
  dimension over $R$; see~\ref{xdim}.

\end{rmk}

We now prove some additional properties of the classes $\F_C$,
$\PP_C$, and $\I_C$. When $R$ is commutative and noetherian, part
(a) appears in~\cite{holm:sdmrghd}[(2.14)].

\begin{prop}   \label{closed}
  Let  ${}_SC_R$ be a semidualizing bimodule.
  \begin{enumerate}[\quad\rm(a)]
  \item\label{FCclosed} The class $\F_C(S)$ is closed under
    coproducts, filtered colimits and summands. If $R$ is right
    coherent, then $\F_C(S)$ is also closed under products.
  \item\label{PCclosed} The class $\PP_C(S)$ is closed under
    coproducts and summands.
  \item\label{ICclosed} The class $\I_C(R)$ is closed under products
    and summands. If $S$ is left noetherian, then $\I_C(R)$ is also
    closed under coproducts and filtered colimits.
  \end{enumerate}
\end{prop}

\begin{proof}
  We only prove~\eqref{FCclosed}, as~\eqref{PCclosed}
  and~\eqref{ICclosed} are similar.  To prove that $\F_C(S)$ is closed
  under summands, consider a split exact sequence of $S$-modules
  \begin{displaymath}
     X= 0 \longrightarrow V' \longrightarrow V
     \longrightarrow V'' \longrightarrow 0,
  \end{displaymath}
  with \mbox{$V \in \F_C(S)$}.  Lemma \ref{uw}\eqref{FC} and
  Proposition~\ref{absums}\eqref{absums1} imply \mbox{$V', V'' \in
    \B_C(S)$}.  The complex $\Hom_S(C,X)$ is split exact and contains
  the module $\Hom_S(C,V)$ as a middle term.  By
  Lemma~\ref{uw}\eqref{FC}, this module is $R$-flat. The flat
  $R$-modules are closed under summands, so the modules $\Hom_S(C,V')$
  and $\Hom_S(C,V'')$ are $R$-flat.  An application of
  Lemma~\ref{uw}\eqref{FC} shows that $V', V'' \in \F_C(S)$.

  The class $\F_C(S)$ is closed under finite direct sums, so in order
  to prove that it is closed under arbitrary coproducts, it suffices
  to prove that $\F_C(S)$ is closed under filtered colimits. Let
  $\{N_\lambda\}$ be a filtered direct system of $C$-flat $S$-modules.
  By Lemma~\ref{uw}\eqref{FC} and
  Proposition~\ref{absums}\eqref{absums1} it follows that $ \varinjlim
  N_\lambda \in \B_C(S)$.  Lemma \ref{uw}\eqref{FC} implies that
  $\{\Hom_S(C,N_\lambda)\}$ is a filtered direct system of flat
  $R$-modules. Since the flat $R$-modules are closed under filtered
  colimits, the module
  \begin{displaymath}
    \Hom_S(C,\varinjlim N_\lambda) \cong
     \varinjlim \Hom_S(C,N_\lambda) 
  \end{displaymath}
  is $R$-flat. The above isomorphism comes from the fact that that
  $\Hom_S(C,-)$ commutes with filtered colimits since ${}_SC$ is
  finitely presented. An application of Lemma~\ref{uw}\eqref{FC}
  implies  $\varinjlim N_\lambda \in \F_C(S)$.

  Finally, if $R$ is right coherent then the class of flat (left)
  $R$-modules is closed under products by \cite[(3.2.24)]{enochs:rha}.
  Since $C_R$ is finitely presented, $C\otimes_R-$ commutes with
  products by \cite[(3.2.22)]{enochs:rha}.  Thus arguments similar to
  the ones above show that $\F_C(S)$ is closed under products.
\end{proof}

 The next result shows, in particular, that the classes $\PP_C(S)$,
 $\F_C(S)$, and $\I_C(R)$ are closed under extensions.

\begin{prop}  \label{cprojext}
  Let  ${}_SC_R$ be a semidualizing bimodule. Consider the
  following exact sequences of $S$- and $R$-modules, respectively
  \begin{align*}
    0\longrightarrow W' \longrightarrow W \longrightarrow
    W'' \longrightarrow 0, \\
  0\longrightarrow U' \longrightarrow U \longrightarrow
    U'' \longrightarrow 0.
  \end{align*}
  The following assertions hold.
  \begin{enumerate}[\quad\rm(a)]
\item\label{f} If $\,W',W'' \in \F_C(S)$, then $W
    \in \F_C(S)$.
\item\label{p} If $\,W',W'' \in \PP_C(S)$, then
    $X$ splits and $W \in \PP_C(S)$.
  \item\label{i} If $\,U',U'' \in \I_C(R)$, then
    $Y$ splits and $U \in \I_C(R)$.
  \end{enumerate}
\end{prop}

\begin{proof}
  We prove only~\eqref{p}, as~\eqref{i} is dual and~\eqref{f} requires only minor adjustments.
  If $P'$ is $R$-projective, then
  Lemma~\ref{tensoreval}\eqref{t3} implies  $\ext^1_S(C,C\otimes_RP')=0$, and hence
   $ \ext^1_S(C,W')=0$
   since $W'$ is $C$-projective.  This, together with
   Lemma~\ref{uw}\eqref{PC}, forces $\Hom_S(C,X)$ to be a split
   exact sequence
   of $R$-projectives.  Thus, the left column in the following commutative diagram,
   $C\otimes_R\Hom_S(C,X)$, is a split exact
   sequence of $C$-projective $S$-modules.
  \begin{equation}
    \tag{\text{$*$}}
    \begin{split}
    \xymatrix{
    0 \ar[d] & 0 \ar[d] \\
    C \otimes_R \Hom_S(C,W') \ar[d]
    \ar[r]^-{\nu^{}_{W'}}_-{\cong} & W' \ar[d]\\
    C \otimes_R \Hom_S(C,W) \ar[d]
    \ar[r]^-{\nu^{}_W} & W \ar[d] \\
    C \otimes_R \Hom_S(C,W'') \ar[d]
    \ar[r]^-{\nu^{}_{W''}}_-{\cong} & W'' \ar[d] \\
    0 & 0
    }
    \end{split}
  \end{equation}
  Lemma~\ref{uw}\eqref{PC} implies that $\nu^{}_{W''}$ and $\nu^{}_{W'}$ are
  isomorphisms, and the five lemma forces $\nu^{}_{W}$ to be an isomorphism
  as well.  Thus, $V\in\PP_C(S)$ the right column in $(*)$ is split exact, as desired.
\end{proof}

\begin{cor}
  Let ${}_SC_R$ be a semidualizing bimodule.  The classes
  $\PP_C^f(S)$, $\F_C^f(S)$, and $\I_C^f(R)$ are closed under
  extensions (see notation in \ref{setup}).\qed
\end{cor}

\begin{para}
  \label{pure}
  Recall that a short exact sequence of $S$-modules
  \begin{displaymath}
    X= 0 \longrightarrow N' \longrightarrow N \longrightarrow
    N'' \longrightarrow 0
  \end{displaymath}
  is \emph{pure exact} if $A\otimes_SX$ is exact for all
  $S^\mathrm{op}$-modules $A$, equivalently, if $\Hom_S(B,X)$ is exact
  for all finitely presented $S$-modules $B$. When $X$ is pure exact,
  $N'$ is a \emph{pure submodule} of $N$, and $N''$ a
  \emph{pure quotient} of $N$.  See~\cite[appendix]{jensen:mta} for
  more details.
\end{para}

We now address how $\F_C(S)$ and $\I_C(R)$ behave with respect to
pure submodules and pure quotients. In the commutative noetherian
setting, this is~\cite[(3.9)]{enochs:cpawac}.

\begin{lem}
  \label{lemce}
  Let ${}_SC_R$ be a faithfully semidualizing bimodule.
  \begin{enumerate}[\quad\rm(a)]
  \item\label{FCcelem} The class $\F_C(S)$ is closed under pure submodules and
    pure quotients.
  \item\label{ICcelem} When $S$ is left noetherian,  the class $\I_C(R)$ is
    closed under pure submodules and pure quotients.
  \end{enumerate}
\end{lem}

\begin{proof}
 \eqref{FCcelem}.
 Consider a pure exact sequence,
 $X$, as in ~\ref{pure}, with
 \mbox{$N \in \F_C(S)$}. Since $C$ is finitely presented over $S$, the
 complex $\Hom_S(C,X)$ is an exact sequence of $R$-modules. We
 claim that $\Hom_S(C,X)$ is pure exact. To this end, let $Q$ be a finitely presented $R$-module.  Since
 $C$ is finitely presented over $S$ and \mbox{$C\otimes_R-$} is right
 exact, the $S$-module $C\otimes_R Q$ is finitely presented.
 By Hom-tensor adjointness
 \begin{displaymath}
   \Hom_R(Q,\Hom_S(C,X)) \cong \Hom_S(C\otimes_R Q,X).
 \end{displaymath}
 It remains to note that the latter complex (and hence also the first)
 is exact since $X$ is pure exact, and $C\otimes_R Q$ is
 finitely presented.

 In the pure exact sequence of $R$-modules $\Hom_S(C,X)$, the
 module $\Hom_S(C,N)$ is $R$-flat by Lemma \ref{uw}\eqref{FC}, since
 $N$ is $C$-flat. Since the class of $R$-flat modules is closed under
 pure submodules and pure quotients, $\Hom_S(C,N')$ and
 $\Hom_S(C,N'')$ are also $R$-flat.  Thus, if we can prove
 \begin{displaymath}
   N'\cong C\otimes_R \Hom_S(C,N')
   \ \ \text{ and } \ \
   N''\cong C\otimes_R \Hom_S(C,N'')
 \end{displaymath}
 then $N',N'' \in \F_C(S)$, as desired. Since $\Hom_R(C,X)$ is pure exact,
 there is a commutative diagram with exact columns
 \begin{displaymath}
    \xymatrix{
    0 \ar[d] & 0 \ar[d] \\
    C \otimes_R \Hom_S(C,N') \ar[d]
    \ar[r]^-{\nu^{}_{N'}} & N' \ar[d]\\
    C \otimes_R \Hom_S(C,N) \ar[d]
    \ar[r]^-{\nu^{}_N}_-{\cong} & N \ar[d] \\
    C \otimes_R \Hom_S(C,N'') \ar[d]
    \ar[r]^-{\nu^{}_{N''}} & N'' \ar[d] \\
    0 & 0.
    }
  \end{displaymath}
  Lemma~\ref{uw}\eqref{FC} implies $N \in \B_C(S)$, so $\nu^{}_N$ is an
   isomorphism.  The snake lemma gives that $\nu^{}_{N'}$
  is injective, $\nu^{}_{N''}$ is surjective, and that
 $   \operatorname{Ker}\,\nu^{}_{N''} \cong
 \operatorname{Coker}\,\nu^{}_{N'}.$
 Thus, it suffices to argue that
  $\operatorname{Ker}\,\nu^{}_{N''}=0$. Since $C$ is faithfully semidualizing, it
  is enough to prove that
 $   \Hom_S(C,\operatorname{Ker}\,\nu^{}_{N''})=0.$
  Applying $\Hom_S(C,-)$ to 
  \begin{displaymath}
    0 \longrightarrow \operatorname{Ker}\,\nu^{}_{N''} \longrightarrow
    C \otimes_R \Hom_S(C,N'')
    \stackrel{\nu^{}_{N''}}{\longrightarrow} N''\longrightarrow 0
  \end{displaymath}
  we see it is enough to show that $\Hom_S(C,\nu^{}_{N''})$ is
  injective. We claim  $\Hom_S(C,\nu^{}_{N''})$ is an
  isomorphism. By Observation~\ref{leftrightinverses}, it suffices to
  argue that $\mu^{}_{\Hom_S(C,N'')}$ is an isomorphism.  Since
  $\Hom_S(C,N'')$ is $R$-flat, this follows from Lemma
  \ref{flat}.

   The proof of \eqref{ICcelem} is dual to that of~\eqref{FCcelem}
  --- using that when $S$ is left noetherian the class
  of injective $S$-modules is closed under pure submodules and pure
  quotients.
\end{proof}


We conclude this section with a result on the (pre)covering and
(pre)enveloping properties of the $C$-flats, $C$-projectives, and
$C$-injectives. In the commutative noetherian setting, parts
\eqref{FCce}, \eqref{ICce} below appear in
\cite[(3.5)]{enochs:cpawac}.

\begin{prop}\label{ce}
Let ${}_SC_R$ be a semidualizing bimodule.
  \begin{enumerate}[\quad\rm(a)]
  \item\label{FCce} The class $\F_C(S)$ is covering on the category of
    $S$-modules.
  \item\label{PCce} The class $\PP_C(S)$ is precovering on the category of
    $S$-modules.
  \item\label{ICce} The class $\I_C(R)$ is enveloping on the category of
    $R$-modules.
\item\label{FCceadd} If $R$ is right coherent and $C$ is faithfully
semidualizing, then the class $\F_C(S)$ is preenveloping on the
category of $S$-modules.
  \item\label{ICceadd} If $S$ is left noetherian and $C$ is
  faithfully semidualizing,
  then the class $\I_C(R)$ is covering on the category of $R$-modules.
  \end{enumerate}
\end{prop}

\begin{proof}
  \eqref{FCce} By Bican-El Bashir-Enochs~\cite[(3)]{BicanBashirEnochs}, the
  class of flat $R$-modules is covering. Thus, for any $S$-module
  $N$,
  the $R$-module $\Hom_S(C,N)$ has an $R$-flat cover
  $  \alpha\colon F \stackrel{}{\longrightarrow} \Hom_S(C,N).$
  Define $\beta$ to be the composite homomorphism
  \begin{displaymath}
    C \otimes_R F \xra{C\otimes_R\alpha} C
    \otimes_R\Hom_S(C,N) \xra{\nu^{}_N} N.
  \end{displaymath}
  This is an $\F_C(S)$-cover of $N$:  To prove the precovering
  property,  consider a homomorphism
  $\tau\colon C \otimes_R G \longrightarrow N$ where $G$ is $R$-flat.
  We want  $\psi\colon C\otimes_R G \longrightarrow C\otimes_R
  F$ such that $\tau=\beta\circ\psi$.  The assignment $\nu^{}$ from \ref{aus}(B3) is
  natural. Also,  Lemma~\ref{uw}\eqref{FC} implies
  $C \otimes_R G \in \B_C(R)$.  Observation~\ref{leftrightinverses} then gives
  rise to a commutative diagram
  \begin{equation}
   \tag{$\dagger$}
\begin{split}
    \xymatrix@C=14ex{C \otimes_R G \ar[r]^-{\tau}
      \ar[d]^-{\cong}_-{\nu^{-1}_{(C\otimes_RG)}=C\otimes_R\mu^{}_G} & N \\
      C \otimes_R\Hom_S(C,C \otimes_R G)
      \ar[r]^-{C\otimes_R\Hom_S(C,\tau)} & C
      \otimes_R\Hom_S(C,N)
      \ar[u]_-{\nu^{}_N}.
    }
    \end{split}
  \end{equation}
  Now, since $\alpha$ is an $R$-flat precover and $G$ is $R$-flat
  there exists a homomorphism $\varphi\colon G \longrightarrow F$
  making the following diagram commute
  \begin{equation}
   \tag{$\dagger\dagger$}
    \begin{split}
    \xymatrix{ G \ar[r]^-{\mu^{}_G}_-{\cong} \ar@{-->}[d]_-{\varphi} &
      \Hom_S(C,C \otimes_R G)
      \ar[d]^-{\Hom_S(C,\tau)} \\
      F \ar[r]^-{\alpha} & \Hom_S(C,N).
    }
    \end{split}
  \end{equation}
  Define $\psi=C\otimes_R\varphi$.  The first equality below comes
  from the commutativity of $(\dagger)$
  \begin{align*}
    \tau &= \nu^{}_N \circ (C\otimes_R\Hom_S(C,\tau)) \circ (C\otimes_R\mu^{}_G) \\
    &= \nu^{}_N\circ(C\otimes_R \alpha) \circ (C\otimes_R \varphi) \\
    &=\beta \circ\psi
  \end{align*}
  while the second comes the commutativity of the diagram induced by
  applying the functor $C\otimes_R -$ to the diagram
  ($\dagger\dagger$). The third holds by the definitions of $\beta$
  and $\psi$, and thus $\beta$ is a precover of $N$, as desired.

  To see that $\beta$ is a cover, let $G=F$, $\tau=\beta$, and
  $\beta\circ\psi=\beta$. We show $\psi$ is an automorphism.
It is straightforward to verify that the following diagram is
commutative, and Lemma~\ref{flat} imples $\mu_F^{}$ is an
isomorphism
 \begin{displaymath}
    \xymatrix@C=14ex{F \ar[r]^-{\alpha}
      \ar[d]^-{\cong}_-{\mu^{}_F} & \Hom_S(C,N) \\
      \Hom_S(C,C \otimes_R F)
      \ar[r]^-{\Hom_S(C,C\otimes_R\alpha)} &
      \Hom_S(C,C\otimes_R\Hom_S(C,N)).
      \ar[u]_-{\Hom_S(C,\nu^{}_N)}
    }
 \end{displaymath}
 The equality $\beta\circ\psi=\beta$ implies
 \begin{displaymath}
   \Hom_S(C,\beta)\circ\Hom(C,\psi)=\Hom_R(C,\beta).
 \end{displaymath}
 Using the diagram immediately above, one checks that the following diagram
 is commutative
  \begin{displaymath}
    \xymatrix{ {} & F \ar[d]^-{\alpha}
      \ar@{-->}[dl]_-{\!\!\!\!\!\!\!\!\!\!\!\!\!\!\!\!\!\!
        \mu_F^{-1}\circ\Hom_S(C,\psi)\circ\mu^{}_F} \\
      F \ar[r]_-{\alpha} & \Hom_S(C,N).
    }
  \end{displaymath}
  Since $\alpha$ is an $R$-flat cover, it follows that
  $\mu_F^{-1}\circ\Hom_S(C,\psi)\circ\mu^{}_F$, and hence $\Hom_S(C,\psi)$ must be
  an automorphism. Finally, the commutative diagram
  \begin{displaymath}
    \xymatrix@C=15ex{C \otimes_R F
      \ar[d]^-{\cong}_-{C\otimes_R\mu^{}_F} \ar[r]^-{\psi} & C
      \otimes_R F \ar[d]_-{\cong}^-{C\otimes_R\mu^{}_F} \\
      C\otimes_R\Hom_S(C,C \otimes_R F)
      \ar[r]_-{C\otimes_R\Hom_S(C,\psi)}^-{\cong} &
      C\otimes_R\Hom_S(C,C \otimes_R F)
    }
  \end{displaymath}
  implies that $\psi$ itself must be an automorphism.

 \eqref{PCce} This is similar to the proof of~\eqref{FCce}, using
  the fact that the class of $R$-projective modules is
  precovering.

  \eqref{ICce} This proof is dual to the proof of~\eqref{FCce}, using that the class of $S$-injective modules is
  enveloping by Xu~\cite[(1.2.11)]{xu:fcm} and Eckmann--Schopf~\cite{EckmannSchopf}.

  \eqref{FCceadd} It suffices by~\cite[(3.5)(c)]{rada:rcecm}
 (see also~\cite[(2.6)(ii)]{holm:cep}) to prove that $\F_C(S)$ is
  closed under arbitrary products and pure submodules, which
  follows immediately from Proposition~\ref{closed}\eqref{FCclosed}
  and Lemma~\ref{lemce}\eqref{FCcelem}.

 \eqref{ICceadd} By~\cite[(2.5)]{holm:cep}, it suffices to show
 that $\I_C(R)$ is closed under coproducts and pure quotients.
 This follows from Proposition~\ref{closed}\eqref{ICclosed} and
  Lemma~\ref{lemce}\eqref{ICcelem}.
\end{proof}

\section{Characterizations of $\A_C$ and $\B_C$ and applications}\label{classifyAB}
Before characterizing the modules in the Auslander and Bass classes
in terms of the $C$-flats, $C$-projectives, and $C$-injectives, we
note an immediate consequence of the adjoint isomorphisms
 \begin{eqnarray*}
    &\Hom_S(C\otimes_R X, I) \cong
    \Hom_R(X,\Hom_S(C,I)), \ \text{and}& \\
    &\Hom_S(C\otimes_R P,Y) \cong
    \Hom_R(P,\Hom_S(C,Y)).
  \end{eqnarray*}
\begin{lem}\label{etcC}
Let ${}_SC_R$ be a semidualizing bimodule, let
 $X$ be a complex of $R$-modules, and let $Y$ a complex of
$S$-modules.
  \begin{enumerate}[\quad\rm(a)]
  \item\label{C1} If $C\otimes_R X$ is exact, then
    $\Hom_R(X,\Hom_S(C,I))$ is exact for all injective $S$-modules
    $I$. Conversely, if $I$ is faithfully $S$-injective and
    $\Hom_R(X,\Hom_S(C,I))$ is exact, then $C\otimes_RX$ is exact.
  \item\label{C2}If $\Hom_S(C,Y)$ is exact, then
    $\Hom_S(C\otimes_RP,Y)$ is exact for all projective $R$-modules
    $P$.  Conversely, if $P$ is faithfully $R$-projective and
    $\Hom_S(C\otimes_R P,Y)$ is exact, then $\Hom_S(C,Y)$ is
    exact.\qed
  \end{enumerate}
\end{lem}

 In the commutative noetherian
setting, Theorem~\ref{intro:acrchar} from the Introduction and
Theorem~\ref{bcschar} appear in~\cite[(3.6, 3.7)]{enochs:cpawac};
see also~\cite[(5.5.4, 5.5.5)]{xu:fcm}. We now prove
Theorem~\ref{intro:acrchar} from the introduction.

\

\emph{Proof of Theorem 2}:
  Assume $M \in \A_C(R)$ so that $\tor_{\geqslant 1}^R(C,M)=0$.
 An augmented projective resolution $P^+$ of $M$
 then gives rise to an exact sequence
  \begin{eqnarray*}
    C\otimes_R P^+ = \ \cdots \longrightarrow
    C\otimes_R P_1 \longrightarrow
    C\otimes_R P_0 \longrightarrow
    C\otimes_R M \longrightarrow 0
  \end{eqnarray*}
 By Proposition~\ref{ce}\eqref{ICce}, the class of
  $C$-injective modules is preenveloping. Thus, $M$ admits an
  augmented proper $C$-injective coresolution.  That is, there is a
  complex
  \begin{eqnarray*}
    U^+= \ 0 \longrightarrow M \longrightarrow U^0
    \longrightarrow U^1 \longrightarrow \cdots
  \end{eqnarray*}
  such that $\Hom_R(U^+,W)$ is exact for \mbox{$W \in
    \I_C(R)$}. In particular, if ${}_RI$ is faithfully injective, then
    $\Hom_R(U^+,\Hom_S(C,I))$ is exact. Thus, Lemma \ref{etcC}\eqref{C1} implies
  \begin{eqnarray*}
    C\otimes_R U^+= \ 0 \longrightarrow C\otimes_R M
    \longrightarrow C\otimes_R U^0
    \longrightarrow C\otimes_R U^1 \longrightarrow \cdots
  \end{eqnarray*}
  is exact. Therefore, we prove that $U^+$ is exact.  The
  complex $X$, obtained by splicing together $P^+$ and $U^+$, then has the
  desired properties.  By Lemma~\ref{uw}\eqref{IC}, $C\otimes_R U^i$
  is injective for all $i\geqslant 0$, so
  \mbox{$C\otimes_R U^+$} is an augmented injective resolution of
  \mbox{$C\otimes_R M$}.  Since the modules $U_i$ and $M$ are in
  $\A_C(R)$, there is an isomorphism 
  \begin{displaymath}
    \Hom_S(C,C\otimes_R U^+) \cong U^+.
  \end{displaymath}
  Since $\ext_S^{\geqslant
    1}(C,C\otimes_R M)=0$, it follows that $U^+$
  is exact.

  Conversely, assume there is a complex $X$ satisfying
   properties (a)-(e) of the theorem, where each
  $P_i$ is $R$-flat. The complex $X$ induces exact sequences
  $P^+$ and $U^+$, as depicted above, and property (e) implies that the complexes
   $C\otimes_R P^+$ and $C \otimes_R U^+$ are exact.

  Since $P_i$ is $R$-flat and $C\otimes_R P^+$ is exact, there is an
  equality $\tor_{\geqslant 1}^R(C,M)=0$.  Also, as
  $C\otimes_R U^+$ is exact and the functor
  $\Hom_S(C,-)$ is left exact, the right column in the following
  commutative diagram is exact
  \begin{displaymath}
    \xymatrix{
      0 \ar[d] & 0 \ar[d] \\
      M \ar[r]^-{\mu^{}_M} \ar[d] &
      \Hom_S(C,C\otimes_R M) \ar[d] \\
      U^0 \ar[r]_-{\cong}^-{\mu^{}_{U^0}} \ar[d] &
      \Hom_S(C,C\otimes_R U^0) \ar[d] \\
      U^1 \ar[r]^-{\mu^{}_{U^1}}_-{\cong} &
      \Hom_S(C,C\otimes_R U^1)
    }
  \end{displaymath}
  By Lemma~\ref{uw}\eqref{IC} one has $U^i \in \A_C(R)$ so
  $\mu^{}_{U^0}$ and $\mu^{}_{U^1}$ are isomorphisms.  The five
  lemma implies
  that $\mu^{}_M$ is an isomorphism.  The $\mu^{}_{U^i}$ and $\mu^{}_M$
  fit together to give an isomorphism of complexes
  $\Hom_S(C,C\otimes_R U^+) \cong U^+$.  Since the
  complex $U^+$ is exact, it follows that
  $\Hom_S(C,C\otimes_R U^+)$ is exact. As
  $C\otimes_R U^+$ is an augmented $S$-injective resolution of
  $C\otimes_R M$, there is an equality
  \begin{displaymath}
    \ext_S^{\geqslant 1}(C,C\otimes_R M)=0.
  \end{displaymath}
  Thus,  $M$ belongs to $\A_C(R)$.\qed

\

The next result is proved in a similar manner.

\begin{thm}
  \label{bcschar}
  Let ${}_SC_R$ be a semidualizing bimodule.  An $S$-module $N$
  belongs to $\B_C(S)$ if and only if there exists a complex
  of $S$-modules
  \begin{eqnarray*}
    Y= \ \cdots \longrightarrow W_1 \longrightarrow W_0
    \longrightarrow I^0 \longrightarrow I^1 \longrightarrow
    \cdots
  \end{eqnarray*}
  that satisfies the following conditions
  \begin{enumerate}[\quad\rm(a)]
\item The complex $Y$ is exact;
  \item Each $I^i$ is $S$-injective;
  \item Each $W_i$ is a $C$-projective (or $C$-flat);
  \item There is an isomorphism $N \cong
    \operatorname{Ker}(I^0 \longrightarrow I^1)$; and
  \item The complex $\Hom_S(C,Y)$ is exact.
  \end{enumerate}
  Moreover, if \mbox{$N\in\B_C(S)$} then any complex constructed by
  splicing together an augmented injective coresolution of $N$ and an
  augmented proper $C$-projective resolution of $N$ will satisfy the
  above properties. \qed
\end{thm}

The next two theorems address the behavior of the classes $\A_C(R)$
and $\B_C(S)$ with respect to short exact sequences,
see~\ref{projres} for the terminology.

\begin{thm}
  \label{ABresolving}
  Let ${}_SC_R$ be a semidualizing bimodule.  The classes
  $\A_C(R)$ and $\A_C^{f}(R)$ are projectively resolving, and the class $\B_C(S)$ is
  injectively resolving.
\end{thm}

\begin{proof}
  We prove that $\A_C(R)$ is projectively resolving; the other proof
  is similar.  By Lemma~\ref{flat}, the class $\A_C(R)$ contains the
  $R$-projective modules.  Thus, it suffices to show that, given an exact
  sequence of $R$-modules
  \begin{displaymath}
   X= 0 \longrightarrow M' \longrightarrow M \longrightarrow
    M'' \longrightarrow 0
  \end{displaymath}
  with \mbox{$M'' \in \A_C(R)$}, then \mbox{$M' \in \A_C(R)$}
  if and only if \mbox{$M \in \A_C(R)$}. Since $M'' \in
  \A_C(R)$ we have \mbox{$\tor_1^R(C,M'')=0$}. In particular,
  the complex
  $C\otimes_R X$ is exact.  By
  Lemma~\ref{etcC}\eqref{C1}, the complex $\Hom_R(X,U)$ is exact for
  all $U \in \I_C(R)$.  Moreover, the class $\I_C(R)$ is closed under
  finite direct sums by Proposition~\ref{closed}\eqref{ICclosed} and
  is preenveloping by Proposition~\ref{ce}\eqref{ICce}.  Thus,
  the Horseshoe Lemma for preenveloping
  classes~\cite[(8.2.2)]{enochs:rha} gives a commutative diagram
  \begin{equation}
    \tag{$\dagger$} \label{ladder1}
    \begin{split}
      \xymatrix{ 0 \ar[r]^{} & M' \ar[d]^{} \ar[r]^{} & M
        \ar[d]^{} \ar[r]^{} & M'' \ar[d]^{} \ar[r]^{} & 0\\
        0 \ar[r]^{} & U' \ar[r]^{} & U'\oplus U''
        \ar[r]^{} & U'' \ar[r]{} & 0 }
    \end{split}
  \end{equation}
  with exact rows and where each vertical map gives rise to an augmented proper
  $\I_C(R)$-coresolution. Similarly, the Horseshoe Lemma for
  projective resolutions yields a commutative diagram with exact rows
  \begin{equation}
    \tag{$\dagger\dagger$} \label{ladder2}
    \begin{split}
      \xymatrix{ 0 \ar[r]^{} & P' \ar[d]^{} \ar[r]^{} &
        P'\oplus P''
        \ar[d]^{} \ar[r]^{} & P'' \ar[d]^{} \ar[r]^{} & 0\\
        0 \ar[r]^{} & M' \ar[r]^{} & M \ar[r]^{} & M''
        \ar[r]{} & 0 }
    \end{split}
  \end{equation}
  where each vertical map gives rise to an augmented projective resolution.
  Splicing these diagrams together provides a degreewise split exact
  sequence of complexes
  \begin{equation}
  \tag{\text{$*$}}
  \begin{split}
  \xymatrix{{} & \vdots \ar[d] & \vdots \ar[d] & \vdots \ar[d] & {} \\
    0 \ar[r] & P'_1 \ar[d] \ar[r] & P'_1\oplus P''_1
    \ar[d] \ar[r] & P''_1 \ar[d] \ar[r] & 0 \\
    0 \ar[r] & P'_0 \ar[d] \ar[r] & P'_0 \oplus P''_0
    \ar[d] \ar[r] & P''_0 \ar[d] \ar[r] & 0 \\
    0 \ar[r] & U'_0 \ar[d] \ar[r] & U'_0 \oplus U''_0
    \ar[d] \ar[r] & U''_0 \ar[d] \ar[r] & 0 \\
    0 \ar[r] & U'_1 \ar[d] \ar[r] & U'_1 \oplus U''_1
    \ar[d] \ar[r] & U''_1 \ar[d] \ar[r] & 0 \\
    {} & \vdots & \vdots & \vdots & {} }
  \end{split}
  \end{equation}
  where $P'_i, P''_i$ are projective, $U'_i, U''_i \in
  \I_C(R)$, and where
  \begin{eqnarray*}
   &M' = \operatorname{Coker}(P'_1 \longrightarrow P'_0)
   \ , \
   M'' = \operatorname{Coker}(P''_1 \longrightarrow
   P''_0),& \\
   &M = \operatorname{Coker}(P'_1 \oplus P''_1
   \longrightarrow P'_0 \oplus P''_0).&
 \end{eqnarray*}
 Two of the (nonzero) modules in the complex $X$ are in $\A_C(R)$, and
 so by Theorem~\ref{intro:acrchar}, two of the complexes in $(*)$ are exact.
 The long exact sequence in homology implies that the third complex is
 also exact. Furthermore, if we apply $C\otimes_R-$ to $(*)$, we get
 another degreewise split exact sequence of complexes, and again,
 since two of these complexes are exact (by assumption),  so is the
 third.  Another application of Theorem~\ref{intro:acrchar} completes the
 proof.
\end{proof}

The next result follows from Lemma~\ref{uw} and
Theorem~\ref{ABresolving} by taking appropriate bounded resolutions
and breaking them up into short exact sequences.

\begin{cor}
  \label{inclusions2}
  Let ${}_SC_R$ be a semidualizing bimodule.  The class $\A_C(R)$
  contains the $R$-modules of finite $\I_C$-injective dimension
  and the class
  $\B_C(S)$ contains the $S$-modules of finite $\F_C$-projective dimension and
  finite $\PP_C$- projective dimension.  \qed
\end{cor}

\begin{thm}
  \label{sesGo}
  Let ${}_SC_R$ be a faithfully semidualizing bimodule.  The classes
   $\A_C(R)$ and $\A_C^f(R)$ are closed under cokernels of monomorphisms and
  $\B_C(S)$ is closed under kernels of epimorphisms.
\end{thm}

\begin{proof}
  We only prove the statement for $\A_C(R)$, as the other statements
  are proved similarly. Consider a short exact sequence of $R$-modules
  \begin{displaymath}
    0 \longrightarrow M' \longrightarrow M \longrightarrow
    M'' \longrightarrow 0
  \end{displaymath}
  with $M'$ and $M$ in $\A_C(R)$.  The proof of
  Theorem~\ref{ABresolving} gives the desired conclusion provided
  \mbox{$\tor_1^R(C,M'')=0$}.  To verify this vanishing, note that $\tor_1^R(C,M)=0$
  as \mbox{$M \in\A_C(R)$}. Hence, there is an exact
  sequence 
  \begin{eqnarray*}
    0 \longrightarrow \tor_1^R(C,M'') \longrightarrow
    C\otimes_RM' \longrightarrow C\otimes_RM \longrightarrow
    C\otimes_R M'' \longrightarrow 0.
  \end{eqnarray*}
Using Remark~\ref{torext}, this is an exact sequence of $S$-modules
homomorphisms.  Applying the functor $\Hom_S(C,-)$ to this sequence
  provides the right-hand exact column in the following commutative
  diagram
  \begin{displaymath}
    \xymatrix{ {} & 0 \ar[d] \\
     {} & \Hom_S(C,\tor_1^R(C,M'')) \ar[d] \\
     M' \ar[r]^-{\mu^{}_{M'}}_-{\cong} \ar[d]
      &\Hom_S(C,C\otimes_R M') \ar[d] \\
     M \ar[r]^-{\mu^{}_M}_-{\cong}
      &\Hom_S(C,C\otimes_R M)
    }
  \end{displaymath}
  As $M' \longrightarrow M$ is injective, a diagram chase shows that
 $   \Hom_S(C,\tor_1^R(C,M''))=0.$
  Since $C$ is faithfully semidualizing, $\tor_1^R(C,M'')$, as desired.
\end{proof}

The next result---which is a non-commutative, non-noetherian version
of \cite[(1.2)]{foxby:gmarm}---follows from Theorem~\ref{sesGo} by
taking appropriate bounded resolutions and breaking them up into short
exact sequences.

\begin{cor}  \label{inclusions}
  Let ${}_SC_R$ be a faithfully semidualizing bimodule.  The class
   $\A_C(R)$ (resp., $\A_C^f(R)$) contains the $R$-modules (resp.,
   finite $R$-modules) of finite flat dimension, and the
    class $\B_C(S)$ contains the $S$-modules of
  finite injective dimension.  \qed
\end{cor}

Theorems~\ref{ABresolving} and~\ref{sesGo} immediately give the
following non-commutative, non-noetherian version of
\cite[(1.3)]{foxby:gmarm} and \cite[(5.5.6), (5.5.7)]{xu:fcm}.

\begin{cor}\label{twothree}
Let ${}_SC_R$ be a faithfully semidualizing bimodule. The classes
$\A_C(R)$ and $\B_C(S)$ have the property that if two of three
modules in a short exact sequence are in the class
  then so is the third.\qed
\end{cor}

 We also have the following result related to Proposition
\ref{cprojext}.

\begin{cor}\label{pfi-resolving}
  Let ${}_SC_R$ be a faithfully semidualizing bimodule.  The classes
  $\PP_C(S)$ and $\F_C(S)$ are projectively resolving
  and the class $\I_C(R)$ is injectively resolving.
\end{cor}

\begin{proof}  We only prove the claim for $\PP_C(S)$.
By Proposition~\ref{cprojext}, we only need to argue
that $\PP_C(S)$ is closed under kernels of epimorphisms.
 Consider an exact sequence of $S$-modules
  \begin{displaymath}  
    0 \longrightarrow W' \longrightarrow W \longrightarrow
    W'' \longrightarrow 0
  \end{displaymath}
  with $W,W'' \in \PP_C(S)$.  By Lemma~\ref{uw}\eqref{PC},
 one has $W,W'' \in \B_C(S)$ and so
  Theorem~\ref{sesGo} implies that $W'\in\B_C(S)$.  It follows that
  $\ext_S^1(C,W')=0$.  The same technique as in the proof of
  Proposition~\ref{cprojext} shows that $W'\in\PP_C(S)$.
\end{proof}

\begin{lem}\label{extvan}
  Let ${}_SC_R$ be a semidualizing bimodule, and let $M$ be an
  $R$-module.
 If \mbox{$\,\ext_S^{\geqslant 1}(C,C\otimes_R M)=0$}, then
  \mbox{$\ext_S^{\geqslant 1}(C\otimes_R P,C\otimes_R M)=0$} for all
  projective $R$-modules $P$.
\end{lem}

\begin{proof}
  Let $I$ be an injective resolution of the $S$-module $C\otimes_R M$
  and $P$ a projective $R$-module. The first and fourth isomorphisms
  below are by definition of  $\ext$
  \begin{align*}
    \ext_S^i(C\otimes_RP,C\otimes_R M)
     &\cong
    \operatorname{H}_{-i}\Hom_S(C\otimes_RP,I) \\
     &\cong
    \operatorname{H}_{-i}\Hom_R(P,\Hom_S(C,I)) \\
     &\cong
    \Hom_R(P,\operatorname{H}_{-i}\Hom_S(C,I)) \\
     &\cong
    \Hom_R(P,\ext_S^i(C,C\otimes_R M))
  \end{align*}
  while the second is by Hom-tensor adjointness, and the third is by
  exactness of the functor $\Hom_R(P,-)$. The desired conclusion follows, as
  $\ext_S^{\geqslant 1}(C,C\otimes_R M)=0$.
\end{proof}

The following result can be thought of as a derived version of Foxby
equivalence. For derived Auslander and Bass classes over a commutative
noetherian ring, cf.~Remark \ref{dualizing}, related results can be
found in e.g.~\cite[(4.5)]{christensen:sdctac}.

\begin{thm}\label{hha}
  Let $M$ and $M'$ be $R$-modules, let $N$ and $N'$ be $S$-modules, let
   $\tilde{N}$ be an $S^{\mathrm{op}}$-module, and let $i\geqslant 0$.
  \begin{enumerate}[\quad\rm(a)]
  \item\label{hha1} If $M\in\A_C(R)$ and $\tor^R_{\geqslant
      1}(C,M')=0$ (e.g., if $M'\in\A_C(R)$), then
    \begin{displaymath}
      \ext_R^i(M',M) \cong
      \ext_S^i(C\otimes_R M',C\otimes_R M).
    \end{displaymath}
  \item\label{hha2} If $N\in\B_C(S)$ and $\ext_S^{\geqslant
      1}(C,N')=0$ (e.g., if $N'\in\B_C(S)$), then
    \begin{displaymath}
      \ext_S^i(N,N') \cong
      \ext_R^i(\Hom_S(C,N),\Hom_S(C,N')).
    \end{displaymath}
  \item\label{hha3} If $N \in \B_C(S)$ and $\tor^S_{\geqslant
      1}(\tilde{N},C)=0$, then
    \begin{displaymath}
      \tor_i^S(\tilde{N},N) \cong
      \tor_i^R(\tilde{N}\otimes_S C,\Hom_S(C,N)).
    \end{displaymath}
  \end{enumerate}
  Each isomorphism defined above is a natural isomorphism of abelian
  groups.
\end{thm}

\begin{proof}  \eqref{hha1}
We proceed by induction on $i$.
  For $i=0$, the first isomorphism below holds since $M\in\A_C(R)$
   while the second is Hom-tensor adjointness
  \begin{align*}
    \Hom_R(M',M)
     &\cong
    \Hom_R(M',\Hom_S(C,C\otimes_R M)) \\
     &\cong
    \Hom_S(C\otimes_R M',C\otimes_R M).
  \end{align*}
    Moreover, these
  isomorphisms are natural in $M'$ and $M$. Next assume that $i>0$.  The induction hypothesis implies that for
 $j<i$ there exist isomorphisms
  \begin{displaymath}
    \tag{\text{$\dagger$}}
    \ext_R^j(L',L) \stackrel{\cong}{\longrightarrow}
    \ext_S^j(C\otimes_R L',C\otimes_R L),
  \end{displaymath}
  which are natural for all $R$-modules \mbox{$L \in \A_C(R)$} and
  $L'$ with \mbox{$\tor^R_{\geqslant 1}(C,L')=0$}. Now, consider
  \mbox{$M \in \A_C(R)$} and $M'$ such that \mbox{$\tor^R_{\geqslant
      1}(C,M')=0$}.  There is a projective $R$-module $P'$ which gives
  rise to an exact sequence
  \begin{displaymath}
   X= 0 \longrightarrow K' \longrightarrow P' \longrightarrow
    M' \longrightarrow 0.
  \end{displaymath}
  The equalities
  \mbox{$\tor^R_{\geqslant 1}(C,M')=0=\tor^R_{\geqslant 1}(C,P')$
  }and the
  appropriate long exact sequence imply \mbox{$\tor^R_{\geqslant
      1}(C,K')=0$}.  Thus, we may apply the induction hypothesis to the
  modules \mbox{$L'=K'$} (or \mbox{$L'=P'$}) and \mbox{$L=M$}.  Since
  \mbox{$\tor^R_1(C,M')=0$}, the complex $X$ induces the exact
  sequence of $S$-modules
  \begin{displaymath}
    C\otimes_R X \ = \quad 0 \longrightarrow C\otimes_R K' \longrightarrow
    C\otimes_R P' \longrightarrow C\otimes_R M'
    \longrightarrow 0.
  \end{displaymath}
  The long exact sequence coming from the complexes $\Hom_R(X,M)$ and
  $\Hom_S(C\otimes_R X ,C\otimes_R M)$ give rise to a
  commutative diagram with exact columns
  \begin{displaymath}
    \xymatrix@R=4ex@C=8ex{
      \vdots \ar[d] & \vdots \ar[d] \\
      \ext_R^{i-1}(P',M) \ar[d] \ar[r]^-{\cong} &
      \ext_S^{i-1}(C\otimes_R P',C\otimes_R M)
      \ar[d] \\
      \ext_R^{i-1}(K',M) \ar[d] \ar[r]^-{\cong} &
      \ext_S^{i-1}(C\otimes_R K',C\otimes_R M)
      \ar[d] \\
      \ext_R^i(M',M) \ar[d]  &
      \ext_S^i(C\otimes_R M',C\otimes_R M)
      \ar[d] \\
      \ext_R^i(P',M) \ar@{=}[d] &
      \ext_S^i(C\otimes_RP',C\otimes_R M)
      \ar@{=}[d] \\
      0 & 0.
    }
  \end{displaymath}
  The right zero follows from Lemma~\ref{extvan}, and the two
  isomorphisms come from the induction hypothesis.  Diagram chasing
 provides a unique isomorphism
  \begin{displaymath}
    \tag{\text{$\star$}}
    \ext_R^i(M',M) \stackrel{\cong}{\longrightarrow}
    \ext_S^i(C\otimes_RM',C\otimes_R M)
  \end{displaymath}
  making the induced diagram commutative. It is straightforward to
  verify that $(\star)$ is natural in $M$ and $M'$.
  Parts~\eqref{hha2} and~\eqref{hha3} have similar proofs.
\end{proof}

\section{Auslander and Bass classes over commutative rings}\label{sec:symmetric}

If $R$ is noetherian then the next two results are
in~\cite[(3.2.9)]{christensen:gd} and \cite[(5.8),
(5.9)]{christensen:sdctac}.

\begin{prop} \label{loc}
  Let $R$ be a commutative ring and $C$ a semidualizing $R$-module. Assume
  $M$ and $N$ are $R$-modules and $\p$ is a prime ideal of $R$.
  \begin{enumerate}[\quad\rm(a)]
    \item \label{loc1} If $M \in \A_C(R)$, then $M_\prim \in
    \A_{C_\prim}(R_\prim)$.
    \item \label{loc2} If $N \in \B_C(R)$, then $N_\prim \in
    \B_{C_\prim}(R_\prim)$.
  \end{enumerate}
\end{prop}

\begin{proof}
  We prove only~\eqref{loc1}, as~\eqref{loc2} is similar. By Observation~\ref{sdmobs}\eqref{sdmobs4}, the
  $R_\prim$-module $C_{\prim}$
  is semidualizing.  If $M \in \A_C(R)$, then for $i>0$
  \begin{eqnarray*}
    &\mathrm{Tor}_i^{R_{\prim}}(C_{\prim},M_{\prim}) \cong
    \Tor_i(C,M)_{\prim} = 0,  \\
    &\mathrm{Ext}^i_{R_{\prim}}(C_{\prim},C_{\prim}\otimes_{R_{\prim}}M_{\prim})
    \cong \Ext^i(C,C\otimes_RM)_{\prim} = 0 ;&
  \end{eqnarray*}
  where the second row uses the assumption that $C$ admits a
  degreewise finite $R$-projective resolution.  Furthermore, the
  commutative diagram
  \begin{eqnarray*}
    \xymatrix{M_{\prim} \ar[d]^-{\cong}_-{(\mu^{}_M)_{\prim}}
      \ar[r]^-{\mu^{}_{(M_{\prim})}} &
      \mathrm{Hom}_{R_{\prim}}(C_{\prim},
      C_{\prim}\otimes_{R_{\prim}}M_{\prim})
      \ar[d]^-{\cong} \\
      \Hom_R(C,C\otimes_RM)_{\prim} \ar[r]_-{\cong} &
      \mathrm{Hom}_{R_{\prim}}(C_{\prim},(C\otimes_RM)_{\prim}) }
  \end{eqnarray*}
  shows that
   $ \mu^{}_{(M_{\prim})}$
  is an isomorphism.  Thus, one has $M_\prim \in
  \A_{C_\prim}(R_\prim)$.
\end{proof}

\begin{prop}  \label{abflatinj}
 If $R$ is commutative and $C$ is a semidualizing $R$-module, then the
  following hold for all $R$-modules $M$ and $N$.
  \begin{enumerate}[\quad\rm(a)]
  \item\label{abflatinj1} \mbox{$M \in \A_C(R) \iff \Hom_R(M,I) \in
      \B_C(R)$} for all injective $R$-modules $I$.
  \item\label{abflatinj1comma5} \mbox{$N \in \B_C(R) \iff \Hom_R(N,I) \in
      \A_C(R)$} for all injective $R$-modules $I$.
  \item\label{abflatinj2} $M \in \A_C(R) \iff M\otimes_RF \in \A_C(R)$
    for all flat $R$-modules $F$.
  \item\label{abflatinj3} $N \in \B_C(R) \iff N\otimes_RF \in \B_C(R)$
    for all flat $R$-modules $F$.
  \end{enumerate}
\end{prop}

\begin{proof}We prove \eqref{abflatinj1}, as \eqref{abflatinj1comma5},
  \eqref{abflatinj2}, and \eqref{abflatinj3} are similar. If $I$ is
  injective then Hom-tensor adjointness gives an isomorphism
  \begin{displaymath}
    \Hom_R(\tor_i^R(C,M),I) \cong \ext_R^i(C,\Hom_R(M,I)).
  \end{displaymath}
  Thus, $\tor_i^R(C,M)=0$ if and only if $\ext_R^i(C,\Hom_R(M,I))=0$ for
  all injective modules $I$. The first isomorphism below follows
  from Lemma~\ref{homeval}\eqref{h2}
  \begin{align*}
    \Hom_R(\ext_R^i(C,C\otimes_RM),I)
    &\cong
    \tor_i^R(C,\Hom_R(C\otimes_RM,I)) \\
    &\cong
    \tor_i^R(C,\Hom_R(C,\Hom_R(M,I)))
  \end{align*}
  while the second is by Hom-tensor adjointness. Hence,
  $\ext_R^i(C,C\otimes_RM)=0$ if and only if
  $\tor_i^R(C,\Hom_R(C,\Hom_R(M,I)))=0$ for all injective modules $I$.
  Finally, there is a commutative diagram
  \begin{eqnarray*}
    \xymatrix{C\otimes_R\Hom_R(C\otimes_RM,I)
      \ar[d]^-{\cong}_-{\theta_{C(C\otimes_RM)I}}
      \ar[rrr]_-{\cong}^-{C\otimes_R(\text{adjointness})} & {} & {} &
      C\otimes_R\Hom_R(C,\Hom_R(M,I)) \ar[d]^-{\nu^{}_{\Hom_R(M,I)}} \\
      \Hom_R(\Hom_R(C,C\otimes_RM),I)
      \ar[rrr]_-{\Hom_R(\mu^{}_M,I)} & {} &
      {} & \Hom_R(M,I) \\
 }
  \end{eqnarray*}
  The left most vertical map is an isomorphism by
  Lemma~\ref{homeval}\eqref{h2}.
  It follows that $\mu^{}_M$ is an
  isomorphism if and only if $\nu^{}_{\Hom_R(M,I)}$ is an isomorphism for
  all injective modules $I$. In conclusion, $M \in \A_C(R)$ if and
  only if $\Hom_R(M,I) \in \B_C(R)$ for all injective modules $I$.
\end{proof}

\begin{rmk}
Each of the statements in the above proposition has a third
equivalent condition.  For example, the statements in (a) are both
equivalent to the following: $\Hom_R(M,I)\in\B_C(R)$ for some
faithfully injective module $I$.
\end{rmk}

\section*{Acknowledgments}

 We thank Sean Sather-Wagstaff for his detailed reading of this
manuscript and many thoughtful suggestions.  We thank Hans-Bj{\o}rn
Foxby and S{\o}ren J{\o}ndrup for useful discussions about
Example~\ref{flatalgebras} and Tom Marley for informing us of
Gruson's Theorem.  We also thank  Frank Moore and Lars Winther
Christensen for their helpful comments.

 \providecommand{\bysame}{\leavevmode\hbox to3em{\hrulefill}\thinspace}
 \providecommand{\MR}{\relax\ifhmode\unskip\space\fi MR }
 \providecommand{\MRhref}[2]{%
   \href{http://www.ams.org/mathscinet-getitem?mr=#1}{#2}
 }
 \providecommand{\href}[2]{#2}


\end{document}